\documentclass[12pt,letterpaper]{amsart}
\usepackage{fullpage,url,amssymb,mathdots}
\usepackage{amsmath,mathrsfs}

\newcommand{\C}{\mathbb C}
\newcommand{\R}{\mathbb R}

\newcommand{\A}{\mathbb A}
\newcommand{\p}{\mathfrak p}

\newcommand{\m}{\mathfrak m}

\renewcommand{\o}{\mathfrak o}

\renewcommand{\a}{\mathfrak a}

\newcommand{\GL}{\operatorname{GL}}

\theoremstyle{plain}
\newtheorem{theorem}{Theorem}[section]     
\newtheorem{corollary}[theorem]{Corollary}
\newtheorem{lemma}[theorem]{Lemma}     
\newtheorem{proposition}[theorem]{Proposition}
\newtheorem{definition}[theorem]{Definition}

\numberwithin{equation}{section}
\theoremstyle{remark}
\newtheorem*{remark}{Remark}
\newtheorem*{remarks}{Remarks}
\newenvironment{psmallmatrix}
{\left(\begin{smallmatrix}}
{\end{smallmatrix}\right)}
\begin{document}    
\title{A converse theorem for GL(2) over global fields without root numbers}
\author{Shantanu Agarwal}

\email{shantanuagr123@gmail.com}

\author{M. Krishnamurthy}

\email{mkrish1939@gmail.com}
\dedicatory{To our parents, with much gratitude}

\begin{abstract}
We prove a converse theorem for $\mathrm{GL}(2)$ over
global fields in which the root numbers in the twisted
functional equations need not be prescribed.
Instead, we allow nonzero scalar factors that are
invariant under unramified twists.
We show that these scalars are independent of the
twisting character and obtain a unique isobaric
automorphic representation with the prescribed
unramified local components.
\end{abstract}

\maketitle
\section{Introduction}\label{intro}
Converse theorems characterize automorphic representations
through the analytic properties of their twisted
$L$-functions. In their usual formulation, the functional
equations involve specified global $\epsilon$-factors,
assembled from the corresponding local factors.
In this paper, we allow these functional equations to
hold up to nonzero scalar factors depending on the
twisting character. Our hypothesis is that the scalar factors 
remain unchanged under unramified twists.
We prove that they are then independent of the twisting
character, and deduce automorphy with the prescribed
unramified local components.
This extends the work of Booker~\cite{B} to global fields.
\begin{theorem}\label{mainthm}
Let $F$ be a global field, i.e., a number field or a function field over a finite field. Let $\A_F$ denote the ring of ad\`{e}les of $F$ and let $\pi=\bigotimes_v\pi_v$ be an irreducible,
admissible generic representation of $\GL_2(\A_F)$ with automorphic
central character $\chi_{\pi}$. For every id\`{e}le class character $\chi$ of $F$ 
which is unramified
at all finite places $v$ where $\pi_v$ is ramified, assume that the products 
$\Lambda(s,\pi\times\chi)=\prod_v L(s,\pi_v\times\chi_v)$ and $\Lambda(s,\widetilde{\pi}\times\chi^{-1})=\prod_v L(s,\widetilde{\pi}_v\times\chi^{-1}_v)$
converge absolutely in some right half plane, continue
to entire functions of finite order, and satisfy the functional equation
\begin{equation}\label{funceq}
\Lambda(s,\pi\times\chi)=\epsilon_{\chi}\epsilon(s,\pi\times\chi)
\Lambda(1-s,\widetilde{\pi}\times\chi^{-1}),
\end{equation}
where $\epsilon(s,\pi\times\chi)$ is the product
of the corresponding local $\epsilon$-factors,
defined in \cite{JL} with respect to a fixed
nontrivial additive character
$\psi=\bigotimes_v\psi_v$ of $F\backslash\A_F$,
and $\epsilon_\chi\in\mathbb C^\times$ is independent
of $s$ and satisfies
$\epsilon_{\chi\omega}=\epsilon_\chi$
for every unramified id\`ele class character $\omega$
in the sense of Definition~\ref{def-unr}. Then there exists a unique isobaric automorphic representation $\Pi=\bigotimes_v\Pi_v$
of $\GL_2(\A_F)$ such that $\pi_v\cong\Pi_v$ for all non-archimedean
places $v$ where $\pi_v$ is unramified.
\end{theorem}
\begin{remarks}\hspace{1pt}
\begin{enumerate}
\item 
The existence of twisted functional equations and the
identification of their precise $\epsilon$-factors are
distinct issues in applications of converse theorems.
For example, in the function field setting,
Grothendieck's theory supplies functional equations for
the twists of the $L$-function of an elliptic curve,
while applying Weil's converse theorem requires a
further comparison with its prescribed functional
equations; see \cite[\S3.2]{Ulmer}.
Our theorem replaces the identification of each scalar
multiplier by a compatibility condition under
unramified twists.

\item If we define the equivalence relation $\chi\sim\chi'\leftrightarrow\chi'\chi^{-1}$ is unramified, then we may unambiguously talk about the number $\epsilon_{[\chi]}$ associated with the equivalence class $[\chi]$.

\item
If $\epsilon_\chi=c\ne0$ is independent of $\chi$,
the conclusion follows from a slight modification of
the proof of the converse theorem in
\cite[\S2]{CKPSS2}.
We explain this modification, together with the passage
to the unique isobaric representation with the prescribed
unramified components, in
Lemma~\ref{lemma-constant-converse}.
A variant allowing a fixed nonzero scalar in the
functional equations is also discussed for
$\mathrm{GL}_3$ in \cite[Theorem 13.7]{JPSS79}.
It would be interesting to extend
Theorem~\ref{mainthm} to $\mathrm{GL}_3$ by showing
that invariance under unramified twists again forces
the multipliers to be constant.

\item\label{rmk-normalization} 
The principal step is to show that the multipliers
$\epsilon_{\chi}$ are independent of $\chi$.
We normalize the auxiliary dual function introduced
below so that the corresponding multipliers become
$\epsilon_\chi/\epsilon_{1}$.
Starting with the case $\chi\sim1$, we use Mellin
inversion and finite Fourier analysis to show that
every normalized multiplier equals~$1$.
Thus the original multipliers are all equal to
$\epsilon_{1}$.

\end{enumerate}
\end{remarks}

Our proof draws on Weil's use of character twists
and Mellin inversion \cite[Theorem~4.3.15]{Miy}
and on Booker's approach to recovering root numbers
\cite{B}. Over a general global field, additional
arguments are needed to handle divisor classes
and the possible vanishing of Whittaker coefficients.
These issues are addressed through a sequence of
lemmas developed in the course of the proof.
We now outline the method.

Suppose
\[
f(z)=\sum_{n\geq1}a_n e^{2\pi inz},
\]
where the coefficients $a_n$ have at most polynomial
growth. Thus $f$ is holomorphic on the upper half
plane $\mathbb H$, and the series converges
absolutely and uniformly on compact subsets of
$\mathbb H$.
Fix integers $N\geq1$ and even $k\geq2$.
Consider the completed $L$-function
\[
\Lambda_N(s,f)
=
N^{s/2}(2\pi)^{-s}\Gamma(s)
\sum_{n\geq1}a_n n^{-s}.
\]
Hecke's Mellin inversion
argument shows that $\Lambda_N(s,f)$ extends to an
entire function, bounded in vertical strips, satisfying
\[
\Lambda_N(s,f)=i^k\Lambda_N(k-s,f),
\]
if and only if
\[
f\left(-\frac1{Nz}\right)=N^{k/2}z^k f(z).
\]
Thus a functional equation yields a transformation
law for $f$. Applying the same principle to
character twists of $f$ is the starting point of the
classical converse-theorem method. The twist
by a Dirichlet character $\chi$ is defined as 
\[
f_\chi(z)=\sum_{n\geq1}a_n\chi(n)e^{2\pi inz}.
\]
For a primitive character $\chi$ modulo $q$, this
can be expressed as an average of translates,
\[
f_\chi(z)
=
\frac{1}{\tau(\overline\chi)}
\sum_{a\bmod q}\overline\chi(a)
f\left(z+\frac aq\right),
\qquad
\tau(\overline\chi)
=
\sum_{a\bmod q}\overline\chi(a)e^{2\pi ia/q}.
\]

We use an analogous construction for Whittaker sums
on $\GL_2(\A_F)$. For $\pi=\bigotimes_v\pi_v$ as in
Theorem~\ref{mainthm}, let $\mathfrak n$ be its
conductor. Retain the global additive character
$\psi=\bigotimes_v\psi_v$ fixed above. At each
non-archimedean place, choose ${\sf d}_v\in F_v^\times$
such that the conductor of $\psi_v$ is
${\sf d}_v^{-1}\o_v$. Throughout the paper, let $S_0$ denote the set of
archimedean places of $F$ when $F$ is a number field.
When $F$ is a function field, fix a place $v_\infty$
at which $\pi_{v_\infty}$ is unramified and $\psi_{v_\infty}$
has conductor $\o_{v_\infty}$, and put
$S_0=\{v_\infty\}$.
Take ${\sf d}_v=1$ at the places in $S_0$, and put
\[
{\sf d}=({\sf d}_v)_v,
\qquad
\mathfrak d=\operatorname{div}^{S_0}({\sf d}).
\]
Let $\xi\mapsto W_\xi$ denote the Whittaker model
of $\pi$ with respect to $\psi$.
See Section~\ref{prelim} for details.

 Write
\[
a(y)=\begin{pmatrix}y&0\\0&1\end{pmatrix},
\qquad
u(x)=\begin{pmatrix}1&x\\0&1\end{pmatrix}.
\]
To a pure tensor $\xi=\bigotimes_v\xi_v$, associate
the Whittaker sum
\[
\phi_\xi(g)
=
\sum_{\gamma\in F^\times}W_\xi(a(\gamma)g).
\]
We fix the normalized essential vectors
$\xi_v=\xi_v^0$ for $v\notin S_0$, with
\[
W_{\xi_v^0}(a({\sf d}_v^{-1}))=1,
\]
and allow $\xi_v$ to vary over the smooth vectors
for $v\in S_0$. With the vectors outside $S_0$ fixed,
we write $\phi_{\xi_{S_0}}$ for $\phi_\xi$, where
$\xi_{S_0}=\bigotimes_{v\in S_0}\xi_v$.

We call an effective divisor $\mathfrak m$
\emph{admissible} if its support is disjoint from
$S_0\cup\operatorname{Supp}(\mathfrak n)$. For such a modulus, put
\[
X(\mathfrak m)
=
\left\{
\text{unitary id\`ele class characters }\omega:
\mathfrak f_\omega^{S_0}\mid\mathfrak m
\right\},
\]
where $\mathfrak f_\omega^{S_0}$ denotes the part
of the conductor of $\omega$ supported outside
$S_0$. No restriction is imposed at the places
in $S_0$.

For an admissible modulus
$\mathfrak m=\prod_{v\in S}\mathfrak p_v^{n_v}$,
put ${\sf w}_v=\varpi_v^{-n_v}$ for $v\in S$
and ${\sf w}_v=0$ otherwise. For $\omega\in X(\mathfrak m)$,
define the twist
\[
\phi_{\omega,\mathfrak m}(g)
=
\int_{\prod_{v\in S}\o_v^\times}
\phi_\xi(g\,u({\sf w}u))\,\omega_S(u)\,d^\times u,
\]
where $\omega_S(u)=\prod_{v\in S}\omega_v(u_v)$
and the measure is normalized so that each
$\o_v^\times$ has unit volume.
The local Mellin calculation then introduces the factor
\[
\int_{\o_v^\times}\omega_v(u)\chi_v(u)^{-1}\,d^\times u
=
\begin{cases}
1,&\omega_v|_{\o_v^\times}=\chi_v|_{\o_v^\times},\\
0,&\text{otherwise}.
\end{cases}
\]
Together with the essential-vector calculation
outside $S\cup S_0$, this isolates precisely the
characters $\chi\sim\omega$ (cf. Lemma~\ref{lemma-m3}). 
The local integrals arising from these twists
and their dual counterparts are not directly
related by the usual local functional equation.
For $v\in S$, we introduce normalized polynomials
$D_{\mathfrak m}(s,\pi_v\times\chi_v)$ and establish
a functional equation for them
(Lemma~\ref{lemma-dirfe}).
This supplies the local identities needed to
deduce, by Mellin inversion, identities between
the twisted Whittaker sums from the assumed
global functional equations.

After normalizing the scalar factors by
$\epsilon_{[1]}$, write $\epsilon_{[\omega]}$
for the resulting factors, so that
$\epsilon_{[1]}=1$. Mellin inversion of the functional
equations, followed by summation over
$X(\mathfrak m)/X(1)$, leads to a finite Fourier
transform. More precisely, put
\[
R_{\mathfrak m}
=\prod_{v\in S}\o_v/\mathfrak p_v^{n_v},
\qquad
\mathcal K_{\mathfrak m}=R_{\mathfrak m}^{\times},
\qquad
\Psi_{\mathfrak m}(x)
=\prod_{v\in S}
\psi_v({\sf d}_v^{-1}\varpi_v^{-n_v}x_v).
\]
Restriction to local units identifies
$X(\mathfrak m)/X(1)$ with
$\widehat{\mathcal K}_{\mathfrak m}$.
Writing $\epsilon_\rho$ for the factor corresponding
to $\rho\in\widehat{\mathcal K}_{\mathfrak m}$, define
\[
\widehat\epsilon_{\mathfrak m}(u)
=
\sum_{\rho\in\widehat{\mathcal K}_{\mathfrak m}}
\epsilon_\rho\rho^{-1}(u),
\qquad
S_{\mathfrak m}(x)
=
\frac{1}{|\mathcal K_{\mathfrak m}|}
\sum_{u\in\mathcal K_{\mathfrak m}}
\widehat\epsilon_{\mathfrak m}(u)
\Psi_{\mathfrak m}((1-u)x).
\]
The central assertion is (see Proposition~\ref{prop-S-one})
\[
S_{\mathfrak m}(x)=1
\qquad\forall x\in R_{\mathfrak m}.
\]

Fix an effective divisor $\mathfrak H$, supported
outside $S_0$ and coprime to $\mathfrak n$, and let
${\sf h}$ be its representing id\`ele, with
components $1$ at $S_0$.
We first treat admissible moduli $\mathfrak m$
coprime to $\mathfrak H$ with trivial ray class
modulo $\mathfrak n\mathfrak H$.
Comparing two Whittaker expansions and varying
the smooth vectors at $S_0$ to isolate individual
terms, we obtain
\[
\widetilde c^{S_0}
  (\gamma{\sf n}{\sf h}{\sf f}_{\mathfrak m})
\left(
S_{\mathfrak m}
  ({\sf d}\gamma{\sf n}{\sf h}{\sf f}_{\mathfrak m})
-1
\right)=0
\]
whenever
${\sf d}\gamma{\sf n}{\sf h}{\sf f}_{\mathfrak m}$
is integral outside $S_0$.
Here $\widetilde c^{S_0}$ denotes the product
of the dual Whittaker coefficients outside $S_0$.
Thus the desired identity holds when the
accompanying coefficient is nonzero.

To pass to arbitrary admissible moduli, we first
use the functional equations for an auxiliary
character twist to obtain an effective divisor
$\mathfrak D$ satisfying
\[
[\mathfrak D]=[\mathfrak n\mathfrak d^2]
\quad\text{in }\operatorname{Cl}(\o_F^{S_0}),
\qquad
(\mathfrak D,\mathfrak m)=1,
\qquad
\widetilde c^{S_0}
  (\mathfrak D\mathfrak d^{-1})\ne0.
\]
We then choose an effective divisor $\mathfrak H$
such that
\[
[\mathfrak H]=[\mathfrak d],
\qquad
(\mathfrak H,\mathfrak n\mathfrak m\mathfrak D)=1.
\]
In particular,
\[
[\mathfrak d\mathfrak n\mathfrak H]
=[\mathfrak D].
\]

Given $x\in R_{\mathfrak m}$, the Chinese remainder
theorem and the local Hecke relations then allow
us to construct a multiple $\mathfrak M$ of
$\mathfrak m$ with trivial ray class modulo
$\mathfrak n\mathfrak H$, together with
$\gamma\in F^\times$ such that
\[
\widetilde c^{S_0}
  (\gamma{\sf n}{\sf h}{\sf f}_{\mathfrak M})\ne0.
\]
The construction also ensures that
\[
y={\sf d}\gamma{\sf n}{\sf h}{\sf f}_{\mathfrak m}
\]
is integral outside $S_0$ and represents $x$ in
$R_{\mathfrak m}$.
The comparison identity therefore gives
\[
S_{\mathfrak M}
  ({\sf d}\gamma{\sf n}{\sf h}{\sf f}_{\mathfrak M})
=1.
\]
The compatibility relation
\[
S_{\mathfrak M}
\left({\sf f}_{\mathfrak M}
      {\sf f}_{\mathfrak m}^{-1}y\right)
=
S_{\mathfrak m}(y),
\]
valid when $y$ is integral at every place dividing
$\mathfrak M$, transfers the resulting identity
back to the prescribed residue class $x$.
This proves $S_{\mathfrak m}(x)=1$ for every
$x\in R_{\mathfrak m}$ and every admissible
modulus $\mathfrak m$.

Finally, finite Fourier analysis shows that
$\epsilon_{[\omega]}=1$ for every
$\omega\in X(\mathfrak m)$ and every admissible
modulus $\mathfrak m$.
Undoing the normalization, the original scalar
factors are independent of the twisting character.
The constant-scalar variant of the converse theorem (see Lemma~\ref{lemma-constant-converse})
then yields the required isobaric automorphic
representation.

\subsection*{Acknowledgements}
The function field case of this work formed part of
the first author's Ph.D. thesis \cite{Shan} at the University
of Iowa. The present paper develops that work further
and treats both number fields and function fields. We thank Andrew Booker for many helpful discussions
and for explaining his work~\cite{B}. We are also grateful to Min Lee for her support and interest in this project.
\section{Preliminaries}
\label{prelim}
Let $F$ be a global field, either a number field or a
function field of one variable over its field of
constants $\mathbb F_q$. If $F$ is a number field, let
$\o_F$ denote its ring of integers. For each place $v$
of $F$, let $F_v$ denote the completion of $F$ at $v$.
For a nonarchimedean place $v$, let $\o_v$ denote the
ring of integers of $F_v$, $\p_v$ its maximal ideal,
and $q_v$ the cardinality of $\o_v/\p_v$. Fix a
uniformizer $\varpi_v$, and normalize the absolute
value by $|\!|\varpi_v|\!|_v=q_v^{-1}$.

Let $\mathscr D_F$ denote the divisor group of $F$,
the free abelian group on its nonarchimedean places.
Its elements are formal sums $\sum_{v<\infty}n_vv$,
where $n_v\in\mathbb Z$ and $n_v=0$ for almost all $v$.
We also write such a divisor multiplicatively as
$\prod_{v<\infty}\p_v^{n_v}$, and refer to $\p_v$ as
a prime divisor. A divisor is called positive if
$n_v\geq0$ for every $v$. We write $\a\succeq0$
for a positive divisor $\a$, and denote the monoid
of positive divisors by $\mathscr D_F^+$.
The zero divisor is written $1$ in multiplicative
notation.

When $F$ is a number field, $\mathscr D_F$ is naturally
identified with the group of nonzero fractional ideals
of $\o_F$, and $\mathscr D_F^+$ with the monoid of
nonzero integral ideals. When $F$ is a function field,
write $q_v=q^{d_v}$ and call $d_v$ the degree of $v$.
The degree of $\prod_v\p_v^{n_v}$ is then
$\sum_v n_vd_v$.
Each id\`ele ${\sf x}=({\sf x}_v)_v\in\A_F^\times$
determines a divisor
\[
{\rm div}({\sf x})
=\prod_{v<\infty}\p_v^{{\rm ord}_v({\sf x}_v)}.
\]
Define $N(\p_v)=q_v$ and extend this multiplicatively,
so that
\[
N(\a)=\prod_{v<\infty}q_v^{n_v}
\qquad
\left(\a=\prod_{v<\infty}\p_v^{n_v}\right).
\]

Let $S_\infty$ denote the set of archimedean places
of $F$. If $F$ is a number field, put $S_0=S_\infty$.
If $F$ is a function field, choose a place $v_\infty$
at which the representation $\pi$ of
Theorem~\ref{mainthm} is unramified, and put
$S_0=\{v_\infty\}$. In this case $S_\infty$ is empty. 
For any finite set of places $S$ containing $S_0$, let
\[
\o_F^S
=
\left\{
x\in F:\operatorname{ord}_v(x)\geq0
\text{ for every }v\notin S
\right\}
\]
denote the ring of $S$-integers in $F$.
In particular, $\o_F^{S_0}=\o_F$ in the number-field
case; in the function-field case, its elements may have
poles only at $v_\infty$. Put
$A_S=\o_F^S$, and let $\mathscr D_F^S$ denote the
group of divisors supported outside $S$. For a divisor $\m$, we write $\m^S$ to denote $\m^S=\prod_{v\notin S}\p_v^{{\rm ord}_v(\m)}$.
Define
\[
\operatorname{Cl}(A_S)
=
\mathscr D_F^S
\big/
\left\{\operatorname{div}^S(x):x\in F^\times\right\}.
\]
We denote the class of $\mathfrak a$ by
$[\mathfrak a]$.

For an effective divisor
$\mathfrak b=\prod_{v\notin S}\mathfrak p_v^{b_v}$,
let $\mathscr D_F^S(\mathfrak b)$ be the subgroup
of divisors supported outside
$S\cup\operatorname{Supp}(\mathfrak b)$, and put
\[
P_{\mathfrak b}(A_S)
=
\left\{
\operatorname{div}^S(x):
x\in F^\times,\quad
x-1\in\mathfrak p_v^{b_v}\o_v
\text{ for every }v\mid\mathfrak b
\right\}.
\]
The ray class group modulo $\mathfrak b$ is
\[
\operatorname{Cl}_{\mathfrak b}(A_S)
=
\mathscr D_F^S(\mathfrak b)/P_{\mathfrak b}(A_S).
\]
No sign conditions are imposed at real places.
We denote the ray class of $\mathfrak a$ by
$[\mathfrak a]_{\mathfrak b}$.
Both groups are finite.
When $S=S_0$, we write simply $A=\o_F^{S_0}$.

For each nonarchimedean place $v$ and integer $n\geq0$,
put
\[
U_v^{(n)}
=
\begin{cases}
\o_v^\times,&n=0,\\
1+\p_v^n,&n>0.
\end{cases}
\]
Let $\chi=\bigotimes_v\chi_v$ be a continuous
quasi-character of $F^\times\backslash\A_F^\times$,
which we call an id\`ele class character.
For each nonarchimedean place $v$, let $a(\chi_v)$
be the least integer $n\geq0$ such that $\chi_v$
is trivial on $U_v^{(n)}$. The conductor of $\chi$
is the positive divisor
\[
\mathfrak f_\chi
=\prod_{v<\infty}\p_v^{a(\chi_v)}.
\]
We represent this divisor by the id\`ele
${\sf f}_\chi$ whose component at a nonarchimedean
place $v$ is $\varpi_v^{a(\chi_v)}$ and whose
archimedean components are~$1$.

If $v\nmid\mathfrak f_\chi$, the value
\[
\chi(\p_v):=\chi_v(\varpi_v)
\]
is independent of the choice of uniformizer.
These values extend multiplicatively to a character
of the group of divisors supported away from
$\mathfrak f_\chi$. We also denote this character by
$\chi$, and extend it to a function on
$\mathscr D_F^+$ by setting $\chi(\a)=0$ whenever
$\a$ and $\mathfrak f_\chi$ have a common prime divisor.

\begin{definition}\label{def-unr}
An id\`ele class character $\chi$ is said to be
\emph{unramified} if $\chi_v$ is trivial on
$\o_v^\times$ for every $v\notin S_0$.
\end{definition}

\begin{remark}
Thus ``unramified'' means unramified outside $S_0$.
In the function field case, ramification at
$v_\infty$ is permitted.
\end{remark}

For a positive divisor $\mathfrak f$ supported outside
$S_0$, put
\[
X(\mathfrak f)
=
\left\{
\text{unitary id\`ele class characters }\omega:
\mathfrak f_\omega^{S_0}\mid\mathfrak f
\right\},
\qquad
\mathfrak f_\omega^{S_0}
=
\prod_{v\notin S_0}\p_v^{a(\omega_v)}.
\]
No local restriction is imposed at the places in $S_0$.
In particular, $X(1)$ is the group of unitary
id\`ele class characters unramified outside $S_0$.

We next fix notation for subgroups of $\GL_2$.
Let ${\rm B}={\rm T}{\rm U}$ be the Borel subgroup
of upper triangular matrices, where ${\rm T}$ is
the diagonal torus and ${\rm U}$ is the upper
unipotent subgroup. Let ${\rm B}_1\subset{\rm B}$
be the mirabolic subgroup consisting of matrices
whose lower-right entry is~$1$.

For each place $v$, put $G_v=\GL_2(F_v)$.
For a finite set $S$ of places, write
\[
F_S=\prod_{v\in S}F_v,\qquad
G_S=\prod_{v\in S}G_v,
\]
and let $\A_F^S$ denote the restricted product
$\prod'_{v\notin S}F_v$. Set
\[
G^S=\GL_2(\A_F^S)
=\prod_{v\notin S}'G_v.
\]
Then
\[
\GL_2(\A_F)=G_S\times G^S.
\]
We identify $G_S$ with its image under the embedding
$g_S\mapsto(g_S,1)$. 

For a function $f$ on $G(\A_F)$ and $g_S\in G_S$,
we write $f(g_S)$ for $f((g_S,1))$.
Conversely, if $f_S$ is a function on $G_S$ and
$g\in G(\A_F)$, we write $f_S(g)$ for $f_S(g_S)$,
where $g_S$ is the projection of $g$ to $G_S$.

When $F$ is a number field, write
\[
F_\infty=\prod_{v\mid\infty}F_v,\qquad
G_\infty=\prod_{v\mid\infty}G_v,
\]
and let $K_\infty=\prod_{v\mid\infty}K_v$, where
each $K_v$ is a fixed maximal compact subgroup
of $G_v$. Let $\mathfrak g_\infty$ be the real Lie algebra
of $G_\infty$, and let
$U(\mathfrak g_{\infty,\mathbb C})$ denote the
universal enveloping algebra of its complexification. In the function field case, $G_\infty$
and $K_\infty$ are understood to be trivial.
In particular, $G_{S_0}=G_\infty$ in the number
field case, whereas $G_{S_0}=G_{v_\infty}$ in the
function field case.

For each nonarchimedean place $v$, put
$K_v=\GL_2(\o_v)$. For an integer $m\geq1$, define
\[
K_{1,v}(\p_v^m)
=
\left\{
g\in\GL_2(\o_v):
g\equiv
\begin{pmatrix}
*&*\\0&1
\end{pmatrix}
\pmod{\p_v^m}
\right\},
\]
and
\[
K_{0,v}(\p_v^m)
=
\left\{
g\in\GL_2(\o_v):
g\equiv
\begin{pmatrix}
*&*\\0&*
\end{pmatrix}
\pmod{\p_v^m}
\right\}.
\]
For $m=0$, put
$K_{1,v}(1)=K_{0,v}(1)=K_v$.
For $m\geq1$, the lower-right entry induces
an isomorphism
\[
K_{0,v}(\p_v^m)/K_{1,v}(\p_v^m)
\simeq(\o_v/\p_v^m)^\times.
\]
In particular, $K_{1,v}(\p_v^m)$ is normal in
$K_{0,v}(\p_v^m)$.

Let $K_f=\prod_{v<\infty}K_v$. For
$\a=\prod_{v<\infty}\p_v^{m_v}\in\mathscr D_F^+$,
put
\[
K_i(\a)
=\prod_{v<\infty}K_{i,v}(\p_v^{m_v}),
\qquad i=0,1.
\]
These are compact open subgroups of the group of
finite adelic points, with
\[
K_1(\a)\subseteq K_0(\a)\subseteq K_f.
\]
We also write
\[
K_i^{S_0}(\a)
=\prod_{v\notin S_0}K_{i,v}(\p_v^{m_v}),
\qquad i=0,1.
\]

\subsection{Whittaker models}
Throughout this paper, we fix a non-trivial additive character $\psi=\otimes\psi_v$ of $F\backslash {\A_F}$ such that $\psi_w(x)=e^{-2\pi ix}$ when $F_w=\R$, $\psi_w(x)=e^{-2\pi i(x+\bar{x})}$ when $F_w=\C$. For each nonarchimedean place $v$,
choose $d_v\in F_v^\times$ such that
\[
\operatorname{cond}(\psi_v)=d_v^{-1}\o_v,
\]
where the conductor is the largest fractional ideal
on which $\psi_v$ is trivial. Let $\delta_v={\rm ord}_v(d_v)$. We may take $d_v=1$
at almost all places. Put ${\sf d}=(d_v)_v$, taking
$d_v=1$ at the archimedean places.
For a nonarchimedean place $v$, the character
\[
\psi_v^\circ(x)=\psi_v(d_v^{-1}x)
\]
has conductor $\o_v$. We use $\psi_v^\circ$ for
local calculations when convenient, while retaining
$\psi$ in all global calculations.

Let $\pi=\bigotimes_v\pi_v$ be an irreducible
admissible representation of $\GL_2(\A_F)$.
At an archimedean place $v$, we work with the smooth
Fr\'echet globalization of moderate growth
$V_{\pi_v}^{\infty}$ of the underlying Harish-Chandra
module $V_{\pi_v}$; see \cite{Be-K}.
When discussing smooth vectors and the action of $G_v$,
we use $\pi_v$ also to denote this globalization.

We say that $\pi$ is \emph{generic} if each $\pi_v$
admits a nonzero Whittaker functional with respect
to $\psi_v$, required to be continuous at an
archimedean place. Equivalently, there is a nonzero
$G_v$-intertwining map
\[
V_{\pi_v}\longrightarrow
{\rm Ind}_{{\rm U}(F_v)}^{G_v}\psi_v,
\qquad
\xi'_v\longmapsto W_{\xi'_v},
\]
where the domain is understood to be
$V_{\pi_v}^{\infty}$ at an archimedean place.
Here the induced representation consists of smooth
functions satisfying
\[
W\left(
\begin{pmatrix}1&x\\0&1\end{pmatrix}g
\right)
=\psi_v(x)W(g),
\]
with $G_v$ acting by right translation.
The intertwining map is unique up to scalar
multiplication, and its image
$\mathcal W(\pi_v,\psi_v)$ is called the
Whittaker model of $\pi_v$ with respect to $\psi_v$. If $\psi'_v$ is another nontrivial additive character
of $F_v$, there is a unique $a\in F_v^\times$ such
that $\psi'_v(x)=\psi_v(a^{-1}x)$. For
$W\in\mathcal W(\pi_v,\psi_v)$, put
\[
(L(a)W)(g)
=
W\left(
\begin{pmatrix}a^{-1}&0\\0&1\end{pmatrix}g
\right).
\]
Then $L(a)$ gives a $G_v$-equivariant bijection
$
\mathcal W(\pi_v,\psi_v)
\longrightarrow
\mathcal W(\pi_v,\psi'_v).
$

We write
\[
\mathcal W(\pi,\psi)
=
\bigotimes_{v}\mathcal W(\pi_v,\psi_v)
\]
for the associated global Whittaker model.
The restricted tensor product is taken with respect
to the spherical Whittaker functions normalized
to have value $1$ at the identity at almost all
nonarchimedean places.
For a pure tensor $\xi=\bigotimes_v\xi_v$, the
corresponding Whittaker function is
\[
W_\xi(g)=\prod_v W_{\xi_v}(g_v).
\]
Similarly, we form
$\mathcal W(\widetilde\pi,\psi)$ for the
contragredient representation
$\widetilde\pi=\bigotimes_v\widetilde\pi_v$.

Fix
\[
w=\begin{pmatrix}0&1\\-1&0\end{pmatrix}.
\]
For $W\in\mathcal W(\pi_v,\psi_v)$, define
\[
\widetilde W(g)=W(w\,{}^t g^{-1}).
\]
Then $W\mapsto\widetilde W$ is a linear bijection
\[
\mathcal W(\pi_v,\psi_v)
\longrightarrow
\mathcal W(\widetilde\pi_v,\psi_v).
\]
With the convention $(R(h)W)(g)=W(gh)$, we have
\[
R(h)\widetilde W_\xi
=
\widetilde W_{\xi'},
\qquad
\xi'=\pi_v({}^t h^{-1})\xi.
\]
Moreover,
\[
\widetilde{L(a)W}
=
\chi_{\pi_v}(a^{-1})L(a)\widetilde W,
\qquad a\in F_v^\times.
\]

\subsection{New vectors and the \textit{essential}
Whittaker function}\label{nv}

Assume that $\pi$ is generic. For each nonarchimedean
place $v$, normalize the multiplicative Haar measure
$d^\times a$ by
$\operatorname{vol}(\o_v^\times,d^\times a)=1$.
Recall that $\psi_v$ is assumed to have conductor $\o_v$.

Let $n(\pi_v)$ be the smallest nonnegative integer
such that $V_{\pi_v}$ contains a nonzero vector fixed
by $K_{1,v}(\p_v^{n(\pi_v)})$. The space of such
vectors is one-dimensional, and its nonzero elements
are called \emph{newvectors}; see \cite{cass}.
The integer $n(\pi_v)$ is the conductor exponent
of $\pi_v$. It also determines the exponent in the
local $\epsilon$-factor \cite{JPSS1}:
\[
\epsilon(s,\pi_v,\psi^{\circ}_v)
=
\epsilon(\pi_v,\psi^{\circ}_v)
q_v^{n(\pi_v)(1/2-s)},
\qquad
\epsilon(\pi_v,\psi^\circ_v)
:=
\epsilon\!\left(\tfrac12,\pi_v,\psi^\circ_v\right).
\]
The conductor of $\pi$ is the positive divisor
\[
\mathfrak n=\prod_{v<\infty}\p_v^{n(\pi_v)}.
\]
Let ${\sf n}=({\sf n}_v)_v$ be the id\`ele defined by
${\sf n}_v=\varpi_v^{n(\pi_v)}$ at nonarchimedean
places and ${\sf n}_v=1$ at archimedean places.

For each nonarchimedean place $v$, choose the
normalized newvector $\xi_v^0$ by requiring
$W_{\xi_v^0}(\begin{psmallmatrix}{\sf d}_v^{-1}&\\&1\end{psmallmatrix})=1$. We call $W_{\xi_v^0}$ the
\emph{essential Whittaker function} and write
\[
W_{\pi_v}=W_{\xi_v^0}
\]
to emphasize the underlying representation.
We use the analogous notation $W_{\widetilde\pi_v}$
for the essential Whittaker function of the
contragredient representation.

For an unramified quasi-character $\chi_v$ of
$F_v^\times$, one has
\begin{equation}\label{eqn-newvector1}
\int_{F_v^\times}
W_{\pi_v}\left(
\begin{psmallmatrix}a{\sf d}_v^{-1}&\\&1\end{psmallmatrix}
\right)
\chi_v(a)|\!|a|\!|_v^{s-1/2}\,d^\times a
=
L(s,\pi_v\otimes\chi_v).
\end{equation}
The integral converges for $\Re(s)$ sufficiently
large, and the identity holds elsewhere by
meromorphic continuation. If $\chi_v$ is ramified,
the integral vanishes, since
$W_{\pi_v}(\begin{psmallmatrix}a&\\&1\end{psmallmatrix})$
is invariant under multiplication of $a$ by
$\o_v^\times$.

Suppose now that $\chi_v$ has conductor exponent
$m_v=a(\chi_v)>0$. Then
\begin{equation}\label{eqn-newvector2}
\int_{F_v^\times}
W_{\pi_v}\left(
\begin{psmallmatrix}
a{\sf d}_v^{-1}&a{\sf d}_v^{-1}\varpi_v^{-m_v}\\&1
\end{psmallmatrix}
\right)
\chi_v(a)|\!|a|\!|_v^{s-1/2}\,d^\times a
=
\mathscr G(\chi_v,\psi_v,\varpi_v^{-m_v-\delta_v}),
\end{equation}
as in \cite{BLK1}, where
\[
\mathscr G(\chi_v,\psi_v,y)
=
\int_{\o_v^\times}\chi_v(u)\psi_v(yu)\,d^\times u.
\]
As above, the integral identity is initially
understood in a right half-plane and then by
meromorphic continuation.

For $m_v>0$, the Gauss sums satisfy
\begin{equation}\label{eqn-gaussrel}
\mathscr G(\chi_v,\psi_v,\varpi_v^{-m_v-\delta_v})
\mathscr G(\chi_v^{-1},\psi_v,\varpi_v^{-m_v-\delta_v})
=
q_v^{-m_v}(1-q_v^{-1})^{-2}\chi_v(-1).
\end{equation}
Put
\[
\kappa_v
=
q_v^{m_v/2}(1-q_v^{-1})
\mathscr G(\chi_v^{-1},\psi_v,\varpi_v^{-m_v-\delta_v}).
\]
Then $|\kappa_v|=1$, and \eqref{eqn-gaussrel} gives
\[
\mathscr G(\chi_v,\psi_v,\varpi_v^{-m_v-\delta_v})
=
\kappa_v^{-1}(1-q_v^{-1})^{-1}
q_v^{-m_v/2}\chi_v(-1).
\]
The relation with the abelian local
$\epsilon$-factor is given by \cite{T}:
\begin{equation}\label{eqn-tate}
\epsilon(s,\chi_v,\psi^\circ_{v})
=
\chi_v(\varpi_v^{m_v})
q_v^{-m_v(s-1/2)}
q_v^{m_v/2}(1-q_v^{-1})
\mathscr G(\chi_v^{-1},\psi_v,\varpi_v^{-m_v-\delta_v}).
\end{equation}

\subsection{The local functional equation}
Fix a place $v$ of $F$ and let $W_v\in {\mathcal W}(\pi_v,\psi_v)$. For a smooth character $\chi_v$ of $F_v^{\times}$ consider the integral
\[
\Psi_v(s; W_v,\chi_v)=\int\limits_{F_v^{\times}}W_v\left(\begin{psmallmatrix}x&\\&1\end{psmallmatrix}\right)\chi_v(x)|\!|x|\!|_v^{s-1/2}
d^{\times}x
\]
originally defined for $Re(s)\gg 0$. 
It follows from \cite{JL} that 
\[
\frac{\Psi_v(s;W_v,\chi_v)}{L(s,\pi_v\times\chi_v)}=\Phi_v(s;W_v,\chi_v)
\]
is an entire function of $s$ for all $v$. If $v$ is non-archimedean, then $\Phi_v(s;W_v,\chi_v)\in {\mathbb C}[q_v^s,q_v^{-s}]$ and for almost all $v$, $\Phi_v(s;W_v,\chi_v)\equiv 1$. We also have the analogous object $\Phi_v(s;\widetilde{W}_v,\chi_v^{-1})$ associated to the pair $(\widetilde{\pi}_v,\chi_v^{-1})$ of contragredient representations. The local functional equation then asserts that (see loc.cit.) 
\begin{equation}\label{lfe}
\Phi_v(1-s;\widetilde{W}_v,\chi^{-1}_v)=\epsilon(s,\pi_v\times\chi_v,\psi_v)\Phi_v(s;W_v,\chi_v)
\end{equation}
for all $W_v$. 
\section{Lemmas}
\subsection{Initial setup}\label{CPS}

Let $\pi=\bigotimes_v\pi_v$ be as in
Theorem~\ref{mainthm}. We first define the auxiliary
functions.
The normalization described in
Remark~(\ref{rmk-normalization}) will be made after
Lemma~\ref{lemma-mi1}.

Let $\Phi(\pi)$ denote the space of functions of
the form
\[
\varphi_{\xi'}(g)
=
\sum_{\gamma\in F^\times}
W_{\xi'}\left(
\begin{psmallmatrix}\gamma&\\&1\end{psmallmatrix}g
\right),
\qquad \xi'\in V_\pi.
\]
The group $\GL_2(\A_F)$ acts on $\Phi(\pi)$ by
right translation, and the map
\[
\xi'\longmapsto\varphi_{\xi'}
\]
is an isomorphism of $\GL_2(\A_F)$-representations.

Retain the set $S_0$ fixed in \S\ref{prelim}.
For each $v\notin S_0$, fix the normalized essential
vector $\xi_v^0\in V_{\pi_v}$ as in \S\ref{nv}.
For $v\in S_0$, let $\xi_v$ be an arbitrary smooth
vector, and write
\[
\xi=\xi^{S_0}\otimes\xi_{S_0},
\qquad
\xi^{S_0}=\bigotimes_{v\notin S_0}\xi_v^0,
\qquad
\xi_{S_0}=\bigotimes_{v\in S_0}\xi_v.
\]
At archimedean places, smooth vectors are understood
to belong to the globalization fixed above.
Consider the function
$\phi=\phi_{\xi_{S_0}}:=\varphi_\xi$.
It satisfies:
\begin{itemize}
\item[(i)]
$\phi(\gamma g)=\phi(g)$ for
$\gamma\in{\rm B}_1(F)$;
\item[(ii)]
$\phi(gz)=\chi_\pi(z)\phi(g)$ for a central element
$z$, identified with an id\`ele in $\A_F^\times$;
\item[(iii)]
$\phi(gk)=\phi(g)$ for
$k\in K_1^{S_0}(\mathfrak n)$.
\end{itemize}

For each $v$, put
\[
\tilde\xi_v
=
\pi_v\left(
\begin{psmallmatrix}{\sf n}_v&\\&1\end{psmallmatrix}
\right)\xi_v,
\]
where $\xi_v=\xi_v^0$ for $v\notin S_0$, and set
$\tilde\xi=\bigotimes_v\tilde\xi_v$.
Define
\[
\widetilde\varphi_\xi(g)
=
\sum_{\gamma\in F^\times}
\widetilde W_{\tilde\xi}\left(
\begin{psmallmatrix}\gamma&\\&1\end{psmallmatrix}g
\right).
\]
One checks that
\[
\widetilde\varphi_\xi(g)
=
\chi_\pi(\det g)^{-1}
\varphi_\xi\left(
g\begin{psmallmatrix}0&1\\-{\sf n}&0\end{psmallmatrix}
\right).
\]
The function
$\widetilde\phi=\widetilde\phi_{\xi_{S_0}}
:=\widetilde\varphi_\xi$
satisfies (i) and (iii), and satisfies (ii) with
$\chi_\pi$ replaced by $\chi_\pi^{-1}$. The function $\phi$ has the Fourier
expansion
\begin{equation}\label{FE}
\phi_{\xi_{S_0}}\left(
\begin{psmallmatrix}y&x\\&1\end{psmallmatrix}
\right)
=
\sum_{\gamma\in F^\times}
c^{S_0}(\gamma y)
W_{\xi_{S_0}}\left(
\begin{psmallmatrix}\gamma y&\\&1\end{psmallmatrix}
\right)
\psi(\gamma x),
\qquad y\in\A_F^\times,\quad x\in\A_F,
\end{equation}
where
\[
c^{S_0}(z)
=
W_{\xi^{S_0}}\left(
\begin{psmallmatrix}z&\\&1\end{psmallmatrix}
\right).
\]
Similarly,
\begin{equation}\label{dual-FE}
\widetilde\phi_{\xi_{S_0}}\left(
\begin{psmallmatrix}y&x\\&1\end{psmallmatrix}
\right)
=
\sum_{\gamma\in F^\times}
\tilde c^{S_0}(\gamma y)
\widetilde W_{\xi_{S_0}}\left(
\begin{psmallmatrix}\gamma y&\\&1\end{psmallmatrix}
\right)
\psi(\gamma x),
\end{equation}
where
\[
\tilde c^{S_0}(z)
=
\widetilde W_{\tilde\xi^{S_0}}\left(
\begin{psmallmatrix}z&\\&1\end{psmallmatrix}
\right).
\]
Here $\tilde\xi_v=\xi_v$ for $v\in S_0$, since
${\sf n}_v=1$ at these places.
For $v\notin S_0$,
$\widetilde W_{\tilde\xi_v}$ is a nonzero scalar
multiple of the essential Whittaker function
$W_{\widetilde\pi_v}$.
Consequently, both $c^{S_0}(z)$ and
$\tilde c^{S_0}(z)$ depend only on
$\m=\operatorname{div}^{S_0}(z)$.
They vanish unless
\[
\operatorname{ord}_v(d_vz_v)\geq0
\qquad\text{for every }v\notin S_0,
\]
or, equivalently, unless
$\operatorname{div}^{S_0}({\sf d}z)$ is effective.
We regard these coefficients as functions on
divisors supported outside $S_0$, and write
$c^{S_0}(\m)$ and $\tilde c^{S_0}(\m)$.

We call
$(\phi,\widetilde\phi)
=(\phi_{\xi_{S_0}},\widetilde\phi_{\xi_{S_0}})$
a \emph{dual pair}.
Their Mellin transforms recover the two sides of
the functional equations for twists by id\`ele
class characters unramified outside $S_0$.

For a local Whittaker function $W$ and a
quasi-character $\chi_v$ of $F_v^\times$, write
\[
\Psi(s;W,\chi_v)
=
\int_{F_v^\times}
W\left(
\begin{psmallmatrix}a&\\&1\end{psmallmatrix}
\right)
\chi_v(a)|\!|a|\!|_v^{s-1/2}\,d^\times a,
\]
initially for $\Re(s)$ sufficiently large and
elsewhere by meromorphic continuation.

For an id\`ele class character $\chi$, put
\[
C_{\psi}^{S_0}(s,\chi)
=
\prod_{v\notin S_0}
\chi_v({\sf d}_v)^{-1}|{\sf d}_v|_v^{1/2-s}.
\]
This is a finite product and non-vanishing for every $s$. For $\Re(s)$
sufficiently large, unfolding and
\eqref{eqn-newvector1} gives
\begin{equation}\label{eqn-Mellin1}
\int_{F^\times\backslash\A_F^\times}
\phi\left(
\begin{psmallmatrix}t&\\&1\end{psmallmatrix}
\right)
\chi(t)|\!|t|\!|^{s-1/2}\,d^\times t
=
C_{\psi}^{S_0}(s,\chi)
\left(
\prod_{v\in S_0}\Psi(s;W_{\xi_v},\chi_v)
\right)
L^{S_0}(s,\pi\times\chi).
\end{equation}
where
\[
L^{S_0}(s,\pi\times\chi)
=
\prod_{v\notin S_0}L(s,\pi_v\times\chi_v).
\]

On the other hand, the definition of $\tilde\xi$
gives
\[
\widetilde W_{\tilde\xi}
\left(
g\begin{psmallmatrix}{\sf n}&\\&1\end{psmallmatrix}
\right)
=
\widetilde W_\xi(g).
\]
In particular,
\[
\widetilde W_{\tilde\xi}
\left(
\begin{psmallmatrix}{\sf n}t^{-1}&\\&1\end{psmallmatrix}
\right)
=
\widetilde W_\xi
\left(
\begin{psmallmatrix}t^{-1}&\\&1\end{psmallmatrix}
\right).
\]
Thus, for $\Re(s)$ sufficiently negative, unfolding
and then changing variables gives
\begin{align}
&\int_{F^\times\backslash\A_F^\times}
\widetilde\phi\left(
\begin{psmallmatrix}{\sf n}t^{-1}&\\&1\end{psmallmatrix}
\right)
\chi(t)|\!|t|\!|^{s-1/2}\,d^\times t
\notag\\
&\qquad=
\int_{\A_F^\times}
\widetilde W_\xi\left(
\begin{psmallmatrix}t^{-1}&\\&1\end{psmallmatrix}
\right)
\chi(t)|\!|t|\!|^{s-1/2}\,d^\times t
\notag\\
&\qquad=
\prod_v
\Psi(1-s;\widetilde W_{\xi_v},\chi_v^{-1}).
\label{eqn-Mellin2}
\end{align}
For $v\notin S_0$, the local functional equation
\eqref{lfe} gives
\[
\begin{aligned}
\Psi(1-s;\widetilde W_{\xi_v},\chi_v^{-1})
&=
\chi_v(d_v)^{-1}|d_v|_v^{1/2-s}
\epsilon(s,\pi_v\times\chi_v,\psi_v)\\
&\qquad\times
L(1-s,\widetilde\pi_v\times\chi_v^{-1}).
\end{aligned}
\]
Put
\[
\epsilon^{S_0}(s,\pi\times\chi,\psi)
=
\prod_{v\notin S_0}
\epsilon(s,\pi_v\times\chi_v,\psi_v).
\]
We conclude that
\begin{align}
&\int_{F^\times\backslash\A_F^\times}
\widetilde\phi\left(
\begin{psmallmatrix}{\sf n}t^{-1}&\\&1\end{psmallmatrix}
\right)
\chi(t)|\!|t|\!|^{s-1/2}\,d^\times t
\notag\\
&\qquad=
C_{\psi}^{S_0}(s,\chi)
\left(
\prod_{v\in S_0}
\Psi(1-s;\widetilde W_{\xi_v},\chi_v^{-1})
\right)
\epsilon^{S_0}(s,\pi\times\chi,\psi)
L^{S_0}(1-s,\widetilde\pi\times\chi^{-1}).
\label{eqn-dualphi}
\end{align}
The Mellin identities are understood beyond their
initial domains of convergence by meromorphic
continuation.

\subsection{The lemmas}
We begin with the following result on Mellin inversion
(cf. \cite[Lemma~11.3.1]{JL}).
For a continuous function $f$ on
$F^\times\backslash\A_F^\times$ and an id\`ele class
character $\chi$, define
\[
\mathscr M(s,f,\chi)
=
\int_{F^\times\backslash\A_F^\times}
f(a)\chi(a)|\!|a|\!|^{s-1/2}\,d^\times a
\]
whenever the integral converges absolutely.

\begin{lemma}\label{lemma-JL}
Let $f_1$ and $f_2$ be continuous functions on
$F^\times\backslash\A_F^\times$.
Suppose there is a constant $c>0$ such that, for
every unitary id\`ele class character $\chi$,
$\mathscr M(s,f_1,\chi)$ converges absolutely for
$\Re(s)>c$, and $\mathscr M(s,f_2,\chi)$ converges
absolutely for $\Re(s)<-c$.
For each $\chi$, assume that the functions represented
by these integrals continue to the same entire
function of $s$, bounded in every vertical strip
of finite width. Then $f_1=f_2$.
\end{lemma}

Following \cite[Chapter~IV]{Weil}, a pair
$({\sf f},{\sf e})\in\A_F^\times\times\A_F$ is
called \emph{reduced} if it satisfies the following
conditions. At each nonarchimedean place $v$,
\[
{\sf f}_v=\varpi_v^{n_v},\qquad n_v\geq0,
\qquad {\sf e}_v\in\o_v,
\qquad
\max\{|\!|{\sf f}_v|\!|_v,|\!|{\sf e}_v|\!|_v\}=1.
\]
In particular, ${\sf e}_v\in\o_v^\times$ whenever
$n_v>0$. At an archimedean place $v$, we require
${\sf f}_v\in\mathbb R_{>0}$ and ${\sf e}_v\in F_v$,
with
\[
{\sf f}_v\overline{{\sf f}_v}
+
{\sf e}_v\overline{{\sf e}_v}=1.
\]
Every $h\in{\rm B}_1(\A_F)$ can be written as
\[
h=
\begin{psmallmatrix}
t{\sf f}&t{\sf e}\\&1
\end{psmallmatrix},
\qquad t\in\A_F^\times,
\]
for a reduced pair $({\sf f},{\sf e})$.
We call
\[
{\rm div}({\sf f})
=
\prod_{v<\infty}\p_v^{n_v}
\]
the \emph{conductor} of the reduced pair.

Given a reduced pair $({\sf f},{\sf e})$, we choose
${\sf e}'\in\A_F$ with
${\sf e}'_v=-\overline{{\sf e}_v}$ at every
archimedean place. Its components at nonarchimedean
places are chosen according to the situation.
If the conductor of $({\sf f},{\sf e})$ is coprime
to $\mathfrak n$, choose ${\sf e}'_v\in\o_v$ such that
\begin{equation}\label{eqn-edisjoint}
|\!|1+{\sf n}_v{\sf e}_v{\sf e}'_v|\!|_v
\leq|\!|{\sf f}_v|\!|_v
\end{equation}
at every nonarchimedean place.
If the conductor of $({\sf f},{\sf e})$ is a
multiple of $\mathfrak n$, choose
${\sf e}'_v\in\o_v$ such that
\begin{equation}\label{eqn-emultiple}
|\!|1+{\sf e}_v{\sf e}'_v|\!|_v
\leq|\!|{\sf f}_v|\!|_v.
\end{equation}

Put
\begin{equation}\label{eqn-feeprime}
k_v=
\begin{psmallmatrix}
(1+{\sf e}_v{\sf e}'_v)/{\sf f}_v&-{\sf e}'_v\\
-{\sf e}_v&{\sf f}_v
\end{psmallmatrix},
\qquad
b_v=
\begin{psmallmatrix}{\sf f}_v&{\sf e}_v\\&1\end{psmallmatrix},
\qquad
b'_v=
\begin{psmallmatrix}{\sf f}_v&{\sf e}'_v\\&1\end{psmallmatrix}.
\end{equation}
In the function field case,
$k_{v_\infty}\in\GL_2(\o_{v_\infty})$.
At an archimedean place, $k_v$ belongs to
$\mathrm{SO}(2)$ if $F_v=\mathbb R$ and to
$\mathrm{SU}(2)$ if $F_v=\mathbb C$.
We take the corresponding maximal compact subgroup
$K_v$ to be $\mathrm O(2)$ or $\mathrm U(2)$,
respectively.

With $w=\begin{psmallmatrix}0&1\\-1&0\end{psmallmatrix}$,
direct calculation gives
\[
b_v=w(b'_vk_v)w^{-1}.
\]
Consequently, writing
\[
k_{S_0}=(k_v)_{v\in S_0},\qquad
b_{S_0}=(b_v)_{v\in S_0},\qquad
b'_{S_0}=(b'_v)_{v\in S_0},
\]
we have
\[
b_{S_0}=w(b'_{S_0}k_{S_0})w^{-1}.
\]
When convenient, we write $k_{S_0}({\sf f},{\sf e})$,
suppressing its dependence on the chosen ${\sf e}'$.




For a reduced pair $({\sf f},{\sf e})$, put
\[
f(t)
=
\phi_{\xi_{S_0}}\left(
\begin{psmallmatrix}
t{\sf f}&t{\sf e}\\&1
\end{psmallmatrix}
\right),
\qquad t\in\A_F^\times.
\]

\begin{lemma}\label{lemma-m1}
Let $({\sf f},{\sf e})$ be a reduced pair whose
conductor $\mathfrak m$ is coprime to $\mathfrak n$.
Then $\mathscr M(s,f,\chi)$ vanishes unless
$\mathfrak f_\chi^{S_0}\mid\mathfrak m^{S_0}$.
When this divisibility holds, it is given by
\[
C_{\psi}^{S_0}(s,\chi)\prod_{\substack{v\notin S_0\\v\nmid\mathfrak m}}
L(s,\pi_v\times\chi_v)
\prod_{\substack{v\notin S_0\\v\mid\mathfrak m}}
\Psi(s;L({\sf d}_v)W_{b_v\cdot\xi_v},\chi_v)
\prod_{v\in S_0}
\Psi(s;W_{b_v\cdot\xi_v},\chi_v).
\]
where
$b_v=\begin{psmallmatrix}{\sf f}_v&{\sf e}_v\\&1\end{psmallmatrix}$.
The formula is initially valid for $\Re(s)$
sufficiently large.
\end{lemma}

\begin{proof}
Unfolding gives
\begin{align*}
\mathscr M(s,f,\chi)
&=
\int_{\A_F^\times}
W_\xi\left(
\begin{psmallmatrix}t{\sf f}&t{\sf e}\\&1\end{psmallmatrix}
\right)
\chi(t)|\!|t|\!|^{s-1/2}\,d^\times t\\
&=
\prod_v\Psi(s;W_{b_v\cdot\xi_v},\chi_v).
\end{align*}
Note that, for any $W$, we have
\begin{equation}\label{eqn-diff}
\Psi(s;W,\chi_v)
=
\chi_v(d_v)^{-1}|d_v|_v^{1/2-s}
\Psi(s;L({\sf d}_v)W,\chi_v)
\end{equation}
If $v\notin S_0$ and $v\nmid\mathfrak m$, then
$b_v\cdot\xi_v=\xi_v$. The local integral vanishes
unless $\chi_v$ is unramified. In the latter case,
changing variables and using \eqref{eqn-diff} and \eqref{eqn-newvector1}
gives
\[
\Psi(s;W_{\xi_v},\chi_v)=
\chi_v(d_v)^{-1}|d_v|_v^{1/2-s}
L(s,\pi_v\times\chi_v).
\]
Suppose $v\notin S_0$ and $v\mid\mathfrak m$.
Write ${\sf f}_v=\varpi_v^{n_v}$, where $n_v>0$.
Since the pair is reduced, ${\sf e}_v\in\o_v^\times$,
and
\[
W_{b_v\cdot\xi_v}\left(
\begin{psmallmatrix}t&\\&1\end{psmallmatrix}
\right)
=
\psi_v(t{\sf e}_v)
W_{\xi_v}\left(
\begin{psmallmatrix}t{\sf f}_v&\\&1\end{psmallmatrix}
\right).
\]
This expression is unchanged when $t$ is replaced
by $tu$, for $u\in1+{\sf f}_v\o_v$. Indeed, the Whittaker function is invariant under multiplication by units. If it is nonzero, then
\[
d_vt{\sf f}_v\in\o_v.
\]
For $u\in1+{\sf f}_v\o_v$, we therefore have
\[
t{\sf e}_v(u-1)\in d_v^{-1}\o_v.
\]
Since $\psi_v$ is trivial on $d_v^{-1}\o_v$,
the additive character is unchanged as well.
Consequently,
\[
\Psi(s;W_{b_v\cdot\xi_v},\chi_v)
=
\chi_v(u)^{-1}
\Psi(s;W_{b_v\cdot\xi_v},\chi_v).
\]
Thus the local integral vanishes unless $\chi_v$
is trivial on $1+\p_v^{n_v}$, equivalently
$a(\chi_v)\leq n_v$. The stated form follows from \eqref{eqn-diff}. 
There is no additional restriction at places
in $S_0$. These observations prove the lemma.
\end{proof}

An analogous result holds for
\[
\widetilde f(t)
=
\chi_\pi({\sf f})
(R(k_{S_0})\widetilde\phi_{\xi_{S_0}})\left(
\begin{psmallmatrix}
{\sf n}t^{-1}{\sf f}&{\sf n}t^{-1}{\sf e}'\\&1
\end{psmallmatrix}
\right).
\]
\begin{lemma}\label{lemma-m2}
With notation as in Lemma~\ref{lemma-m1}, and
${\sf e}'$ chosen as in \eqref{eqn-edisjoint},
the Mellin transform
$\mathscr M(s,\widetilde f,\chi)$ vanishes unless
\[
\mathfrak f_\chi^{S_0}\mid\mathfrak m^{S_0}.
\]
When this divisibility holds, it equals
\begin{align*}
&C_\psi^{S_0}(s,\chi)
\prod_{\substack{v\notin S_0\\v\nmid\mathfrak m}}
\left(
\epsilon(s,\pi_v\times\chi_v,\psi_v)
L(1-s,\widetilde\pi_v\times\chi_v^{-1})
\right)\\
&\quad\times
\prod_{\substack{v\notin S_0\\v\mid\mathfrak m}}
\left(
(\chi_{\pi_v}\chi_v^2)({\sf d}_v)
|{\sf d}_v|_v^{2s-1}
\Psi\!\left(
1-s;
\widetilde{L({\sf d}_v)W_{b_v\cdot\xi_v}},
\chi_v^{-1}
\right)
\right)\\
&\quad\times
\prod_{v\in S_0}
\Psi\!\left(
1-s;\widetilde W_{b_v\cdot\xi_v},\chi_v^{-1}
\right),
\end{align*}
where
\[
b_v=
\begin{psmallmatrix}{\sf f}_v&{\sf e}_v\\&1\end{psmallmatrix}.
\]
The formula is initially valid for $\Re(s)$
sufficiently negative.
\end{lemma}
\begin{proof}
Unfolding and the local matrix calculations used
in the proof of Lemma~\ref{lemma-m4}, with the
averaging operators omitted, give
\[
\mathscr M(s,\widetilde f,\chi)
=
\prod_v
\Psi\!\left(
1-s;\widetilde W_{b_v\cdot\xi_v},\chi_v^{-1}
\right).
\]
Here the factors $\chi_{\pi_v}({\sf f}_v)^{-1}$
arising from those calculations cancel the
prefactor $\chi_\pi({\sf f})$ in the definition
of $\widetilde f$.

For $v\notin S_0$, put
\[
C_v(s,\chi_v)
=
\chi_v({\sf d}_v)^{-1}
|{\sf d}_v|_v^{1/2-s}.
\]
If $v\nmid\mathfrak m$, then
$b_v\cdot\xi_v=\xi_v$.
The local functional equation and
Lemma~\ref{lemma-m1} show that the local integral
vanishes unless $\chi_v$ is unramified.
When $\chi_v$ is unramified, it equals
\[
C_v(s,\chi_v)
\epsilon(s,\pi_v\times\chi_v,\psi_v)
L(1-s,\widetilde\pi_v\times\chi_v^{-1}).
\]

Now suppose $v\notin S_0$ and $v\mid\mathfrak m$.
For any $W\in\mathcal W(\pi_v,\psi_v)$,
a change of variables gives
\[
\Psi(1-s;\widetilde W,\chi_v^{-1})
=
C_v(s,\chi_v)^{-1}
\Psi(1-s;L({\sf d}_v)\widetilde W,\chi_v^{-1}).
\]
Using
\[
L({\sf d}_v)\widetilde W
=
\chi_{\pi_v}({\sf d}_v)
\widetilde{L({\sf d}_v)W},
\]
we obtain
\begin{align*}
\Psi(1-s;\widetilde W,\chi_v^{-1})
&=
C_v(s,\chi_v)
(\chi_{\pi_v}\chi_v^2)({\sf d}_v)
|{\sf d}_v|_v^{2s-1}\\
&\qquad\times
\Psi\!\left(
1-s;\widetilde{L({\sf d}_v)W},\chi_v^{-1}
\right).
\end{align*}
Apply this with $W=W_{b_v\cdot\xi_v}$.
Moreover, the local functional equation transfers
the vanishing established in Lemma~\ref{lemma-m1}
to the dual integral. Thus this factor vanishes
unless
\[
a(\chi_v)\leq\operatorname{ord}_v(\mathfrak m).
\]

At the places in $S_0$, we retain the local
integrals without modification.
Multiplying the local identities and using
\[
\prod_{v\notin S_0}C_v(s,\chi_v)
=
C_\psi^{S_0}(s,\chi)
\]
proves the lemma.
\end{proof}

\begin{lemma}\label{lemma-mi1}
Under the hypotheses of Theorem~\ref{mainthm}, let
$(\phi,\widetilde\phi)$ be associated to any choice
of smooth vectors at the places in $S_0$.
Let $({\sf f},{\sf e})$ be a reduced pair whose
conductor is supported in $S_0$, and choose
${\sf e}'$ satisfying \eqref{eqn-edisjoint}.
Put $\epsilon=\epsilon_{[1]}$.
Then, for every $t\in\A_F^\times$,
\[
\phi\left(
\begin{psmallmatrix}t{\sf f}&t{\sf e}\\&1\end{psmallmatrix}
\right)
=
\epsilon\,\chi_\pi({\sf f})
(R(k_{S_0})\widetilde\phi)\left(
\begin{psmallmatrix}
{\sf n}t^{-1}{\sf f}&{\sf n}t^{-1}{\sf e}'\\&1
\end{psmallmatrix}
\right),
\]
where $k_{S_0}=k_{S_0}({\sf f},{\sf e})$.
\end{lemma}
\begin{proof}
Let $f$ and $\widetilde f$ be as in
Lemmas~\ref{lemma-m1} and~\ref{lemma-m2}.
Since the conductor of $({\sf f},{\sf e})$ is
supported in $S_0$, both Mellin transforms vanish
unless $\chi$ is unramified outside $S_0$.
For such a character, the hypotheses of
Theorem~\ref{mainthm} give
$\epsilon_\chi=\epsilon_{[1]}=\epsilon$.

Put
\[
b_v=
\begin{psmallmatrix}{\sf f}_v&{\sf e}_v\\&1\end{psmallmatrix}.
\]
By Lemma~\ref{lemma-m1},
\[
\mathscr M(s,f,\chi)
=
C_\psi^{S_0}(s,\chi)
\Lambda(s,\pi\times\chi)
\prod_{v\in S_0}
\frac{\Psi(s;W_{b_v\cdot\xi_v},\chi_v)}
     {L(s,\pi_v\times\chi_v)}.
\]
Similarly, Lemma~\ref{lemma-m2} and the local
functional equations at the places in $S_0$ give
\[
\begin{aligned}
\mathscr M(s,\widetilde f,\chi)
&=
C_\psi^{S_0}(s,\chi)
\epsilon(s,\pi\times\chi,\psi)
\Lambda(1-s,\widetilde\pi\times\chi^{-1})\\
&\qquad\times
\prod_{v\in S_0}
\frac{\Psi(s;W_{b_v\cdot\xi_v},\chi_v)}
     {L(s,\pi_v\times\chi_v)}.
\end{aligned}
\]
The assumed global functional equation therefore
implies
\[
\mathscr M(s,f,\chi)
=
\epsilon\,\mathscr M(s,\widetilde f,\chi)
\]
as meromorphically continued functions.
Together with the vanishing when $\chi$ is
ramified outside $S_0$, this establishes the
identity for every unitary id\`ele class
character $\chi$.

Applying Lemma~\ref{lemma-JL} to $f$ and
$\epsilon\widetilde f$ gives
$f(t)=\epsilon\widetilde f(t)$, as required.
\end{proof}

In Lemma~\ref{lemma-mi1}, we used only the functional
equations for twists by characters unramified
outside $S_0$. To extract information from the
remaining functional equations by Mellin inversion,
we introduce twists of the auxiliary functions.

Fix a positive divisor $\m$ supported away from
$\mathfrak n$ and $S_0$, and put $S={\rm Supp}(\m)$. Henceforth, we call such a modulus {\it admissible}. 
Write $n_v={\rm ord}_v(\m)$ and choose the id\`ele
${\sf f}_{\m}$ with components
\[
{\sf f}_{\m,v}
=
\begin{cases}
\varpi_v^{n_v},&v\in S,\\
1,&v\notin S.
\end{cases}
\]
For $\omega\in X(\m)$, define local vectors by
$\xi_{\omega,v}=\xi_v^0$ for $v\notin S\cup S_0$,
by
\[
\xi_{\omega,v}
=
\int_{\o_v^\times}
\omega_v(u)
\pi_v\left(
\begin{psmallmatrix}
1&u\varpi_v^{-n_v}\\&1
\end{psmallmatrix}
\right)
\xi_v^0\,d^\times u
\]
for $v\in S$, and by $\xi_{\omega,v}=\xi_v$ for
$v\in S_0$. Here each $\o_v^\times$ has unit volume,
and $\xi_v$ is the chosen smooth vector at $v\in S_0$.
Put $\xi_\omega=\bigotimes_v\xi_{\omega,v}$ and define
\[
\phi_\omega(g)
=
\phi_{\omega,\m,\xi_{S_0}}(g)
:=
\varphi_{\xi_\omega}(g)
=
\sum_{\gamma\in F^\times}
W_{\xi_\omega}\left(
\begin{psmallmatrix}\gamma&\\&1\end{psmallmatrix}g
\right).
\]

Let ${\sf w}\in\A_F$ have components
${\sf w}_v=\varpi_v^{-n_v}$ for $v\in S$ and
${\sf w}_v=0$ otherwise. Write
\[
\omega_S(u)=\prod_{v\in S}\omega_v(u_v),
\qquad u\in\prod_{v\in S}\o_v^\times.
\]
Then
\begin{equation}\label{eqn-phiglobal}
\phi_\omega(g)
=
\sum_{\gamma\in F^\times}
\int_{\prod_{v\in S}\o_v^\times}
W_\xi\left(
\begin{psmallmatrix}\gamma&\\&1\end{psmallmatrix}
g
\begin{psmallmatrix}1&{\sf w}u\\&1\end{psmallmatrix}
\right)
\omega_S(u)\,d^\times u,
\end{equation}
where ${\sf w}u$ has component ${\sf w}_vu_v$
at $v\in S$ and zero elsewhere.

The construction depends only on
$\omega_v|_{\o_v^\times}$ for $v\in S$.
Consequently, it depends only on the class
$[\omega]\in X(\m)/X(1)$, and we also write
$\phi_{[\omega]}$.
Its Fourier expansion is
\[
\phi_{[\omega]}\left(
\begin{psmallmatrix}{\sf y}&{\sf x}\\&1\end{psmallmatrix}
\right)
=
\sum_{\gamma\in F^\times}
W_\xi\left(
\begin{psmallmatrix}\gamma{\sf y}&\\&1\end{psmallmatrix}
\right)
\left(
\prod_{v\in S}
\mathscr G(\omega_v,\psi_v,
          \gamma{\sf y}_v\varpi_v^{-n_v})
\right)
\psi(\gamma{\sf x}).
\]
Thus the Fourier coefficients of $\phi_{[\omega]}$
are obtained from those of $\phi$ by multiplication
by the indicated product of Gauss integrals.

We define the twist of the normalized dual function
by
\begin{equation}\label{eqn-dualphiglobal}
\widetilde\phi_\omega(g)
=
\int_{\prod_{v\in S}\o_v^\times}
\widetilde\phi\left(
g\begin{psmallmatrix}1&{\sf w}u\\&1\end{psmallmatrix}
\right)
\omega_S(u)\,d^\times u.
\end{equation}
In terms of the Whittaker function,
this is
\[
\widetilde\phi_\omega(g)
=
\sum_{\gamma\in F^\times}
\int_{\prod_{v\in S}\o_v^\times}
\widetilde W_{\tilde\xi}\left(
\begin{psmallmatrix}\gamma&\\&1\end{psmallmatrix}
g
\begin{psmallmatrix}1&{\sf w}u\\&1\end{psmallmatrix}
\right)
\omega_S(u)\,d^\times u.
\]

For $v\in S$, put
\[
\tilde\xi_{\omega,v}
=
\int_{\o_v^\times}
\omega_v(u)
\pi_v\left(
\begin{psmallmatrix}
1&0\\-{\sf w}_vu&1
\end{psmallmatrix}
\right)
\tilde\xi_v\,d^\times u,
\]
and put $\tilde\xi_{\omega,v}=\tilde\xi_v$ for
$v\notin S$.
Then
\[
\widetilde\phi_\omega(g)
=
\sum_{\gamma\in F^\times}
\widetilde W_{\tilde\xi_\omega}\left(
\begin{psmallmatrix}\gamma&\\&1\end{psmallmatrix}g
\right),
\qquad
\tilde\xi_\omega=\bigotimes_v\tilde\xi_{\omega,v}.
\]
This construction also depends only on $[\omega]$.

Finally, for $({\sf f},{\sf e})$ and ${\sf e}'$
as in Lemma~\ref{lemma-mi1}, put
\[
f_{[\omega]}(t)
=
\phi_\omega\left(
\begin{psmallmatrix}
t{\sf f}&t{\sf e}\\&1
\end{psmallmatrix}
\right)
\]
and
\[
\widetilde f_{[\omega]}(t)
=
\chi_\pi(\m)\chi_\pi({\sf f})
(R(k_{S_0})\widetilde\phi_{\omega^{-1}})
\left(
\begin{psmallmatrix}
{\sf n}{\sf f}_{\m}^{2}t^{-1}{\sf f}
&
{\sf n}{\sf f}_{\m}^{2}t^{-1}{\sf e}'\\
&1
\end{psmallmatrix}
\right).
\]

\begin{lemma}\label{lemma-m3}
The Mellin transform ${\mathscr M}(s,f_{[\omega]},\chi)$
vanishes unless $\chi\sim\omega$. When this condition holds,
it equals
\[
C_\psi^{S_0}(s,\chi)
\prod_{v\notin S\cup S_0}
L(s,\pi_v\times\chi_v)
\prod_{v\in S}
\Psi(s;L({\sf d}_v)W_{u_v\cdot\xi_v},\chi_v)
\prod_{v\in S_0}
\Psi(s;W_{b_v\cdot\xi_v},\chi_v).
\]
where
\[
b_v=\begin{psmallmatrix}{\sf f}_v&{\sf e}_v\\&1\end{psmallmatrix},
\qquad
u_v=\begin{psmallmatrix}1&\varpi_v^{-n_v}\\&1\end{psmallmatrix},
\qquad
n_v=\operatorname{ord}_v(\mathfrak m).
\]
\end{lemma}

\begin{proof}
We first work with $\operatorname{Re}(s)$ sufficiently large
that the integrals below converge absolutely.
As in the proof of Lemma~\ref{lemma-m1}, unfolding gives
\[
{\mathscr M}(s,f_{[\omega]},\chi)
=
\int_{F^\times\backslash\A_F^\times}
f_{[\omega]}(t)\chi(t)|\!|t|\!|^{s-1/2}\,d^\times t
=
\prod_v
\Psi\bigl(s;W_{b_v\cdot\xi_{\omega,v}},\chi_v\bigr).
\]

Since the conductor of $({\sf f},{\sf e})$ is supported
in $S_0$, we have ${\sf f}_v=1$ and ${\sf e}_v\in\o_v$
for every $v\notin S_0$.

If $v\notin S\cup S_0$, then
$\xi_{\omega,v}=\xi_v=\xi_v^0$ and
$b_v=u({\sf e}_v)$ fixes $\xi_v$.
Consequently,
\[
\Psi\bigl(s;W_{b_v\cdot\xi_{\omega,v}},\chi_v\bigr)
=
\begin{cases}
\chi_v({\sf d}_v)^{-1}|{\sf d}_v|_v^{1/2-s}
L(s,\pi_v\times\chi_v),
    & \text{if $\chi_v$ is unramified},\\
0,  & \text{otherwise}.
\end{cases}
\]

Now let $v\in S$. The vector $\xi_v=\xi_v^0$ is
spherical, since $S$ is disjoint from the conductor
of $\pi$. Moreover, $b_v=u({\sf e}_v)$ fixes $\xi_v$
and commutes with every $u(\varpi_v^{-n_v}u)$.
It therefore fixes $\xi_{\omega,v}$ as well.
By the definition of $\xi_{\omega,v}$,
\begin{align}\label{eqn-double}
\Psi(s;W_{b_v\cdot\xi_{\omega,v}},\chi_v)
&=
\int_{F_v^\times}\int_{\o_v^\times}
\omega_v(u)\,
W_{\xi_v}\!\left(
\begin{psmallmatrix}t&\\&1\end{psmallmatrix}
\begin{psmallmatrix}1&\varpi_v^{-n_v}u\\&1\end{psmallmatrix}
\right)
\chi_v(t)|\!|t|\!|_v^{s-1/2}
\,d^\times u\,d^\times t
\notag\\
&=
\int_{F_v^\times}\int_{\o_v^\times}
\omega_v(u)\,
W_{\xi_v}\!\left(
\begin{psmallmatrix}tu&\\&1\end{psmallmatrix}
\begin{psmallmatrix}1&\varpi_v^{-n_v}\\&1\end{psmallmatrix}
\right)
\chi_v(t)|\!|t|\!|_v^{s-1/2}
\,d^\times u\,d^\times t.
\end{align}
The second equality follows from
\[
\begin{psmallmatrix}t&\\&1\end{psmallmatrix}
\begin{psmallmatrix}1&\varpi_v^{-n_v}u\\&1\end{psmallmatrix}
=
\begin{psmallmatrix}tu&\\&1\end{psmallmatrix}
\begin{psmallmatrix}1&\varpi_v^{-n_v}\\&1\end{psmallmatrix}
\begin{psmallmatrix}u^{-1}&\\&1\end{psmallmatrix},
\]
because the last matrix fixes $\xi_v$.

Interchanging the integrals and replacing $t$ by $tu^{-1}$
in the inner integral, we obtain
\[
\Psi(s;W_{\xi_{\omega,v}},\chi_v)
=
\left(
\int_{\o_v^\times}
\omega_v(u)\chi_v^{-1}(u)\,d^\times u
\right)
\Psi(s;W_{u_v\cdot\xi_v},\chi_v).
\]
Since $\o_v^\times$ has volume $1$, character
orthogonality gives
\[
\int_{\o_v^\times}
\omega_v(u)\chi_v^{-1}(u)\,d^\times u
=
\begin{cases}
1, & \text{if $\chi_v\omega_v^{-1}$ is unramified},\\
0, & \text{otherwise}.
\end{cases}
\]
In particular,
\begin{equation}\label{eqn-localS}
\Psi(s;W_{\xi_{\omega,v}},\chi_v)
=
\Psi(s;W_{u_v\cdot\xi_v},\chi_v)
\end{equation}
whenever $\chi_v\omega_v^{-1}$ is unramified.
Changing variables 
\[
\Psi(s;W_{u_v\cdot\xi_v},\chi_v)
=
\chi_v({\sf d}_v)^{-1}|{\sf d}_v|_v^{1/2-s}
\Psi(s;L({\sf d}_v)W_{u_v\cdot\xi_v},\chi_v).
\]

Finally, if $v\in S_0$, then $\xi_{\omega,v}=\xi_v$,
so the local factor remains
\[
\Psi(s;W_{b_v\cdot\xi_v},\chi_v).
\]

Since $\omega$ is unramified outside $S\cup S_0$,
the conditions obtained above are together equivalent to
$\chi\omega^{-1}$ being unramified outside $S_0$,
that is, to $\chi\sim\omega$.
Multiplying the local factors 
proves the assertion in the domain of absolute
convergence, and hence wherever the Mellin transform
admits meromorphic continuation.
\end{proof}
Thus the unipotent average defining $\phi_{[\omega]}$
isolates the twists by characters $\chi\sim\omega$.

\begin{lemma}\label{lemma-m4}
The Mellin transform
${\mathscr M}(s,\widetilde f_{[\omega]},\chi)$
vanishes unless $\chi\sim\omega$. When this condition holds,
it equals
\begin{align*}
&C_\psi^{S_0}(s,\chi)
\prod_{v\notin S\cup S_0}
\left(
\epsilon(s,\pi_v\times\chi_v,\psi_v)\,
L(1-s,\widetilde\pi_v\times\chi_v^{-1})
\right)\\
&\quad\times
\prod_{v\in S}
\left(
(\chi_{\pi_v}\chi_v^2)(\varpi_v^{n_v})\,
q_v^{-n_v(2s-1)}
(\chi_{\pi_v}\chi_v^2)({\sf d}_v)\,
|{\sf d}_v|_v^{2s-1}
\Psi\!\left(
1-s;
\widetilde{L({\sf d}_v)W_{\bar u_v\cdot\xi_v}},
\chi_v^{-1}
\right)
\right)\\
&\quad\times
\prod_{v\in S_0}
\Psi\!\left(
1-s;\widetilde W_{b_v\cdot\xi_v},\chi_v^{-1}
\right).
\end{align*}
where
\[
\bar u_v=
\begin{psmallmatrix}
1&\\-\varpi_v^{-n_v}&1
\end{psmallmatrix},
\qquad
n_v=\operatorname{ord}_v(\mathfrak m).
\]
\end{lemma}

\begin{proof}
We first work with $\operatorname{Re}(s)$ sufficiently
negative that the integrals below converge absolutely.
Unfolding gives
\[
{\mathscr M}(s,\widetilde f_{[\omega]},\chi)
=
\chi_\pi(\mathfrak m)\chi_\pi({\sf f})
\prod_v I_v,
\]
where
\[
I_v=
\begin{cases}
\displaystyle
\int_{F_v^\times}
\widetilde W_{\tilde\xi_{\omega^{-1},v}}
\!\left(
\begin{psmallmatrix}
{\sf n}_v{\sf f}_{\mathfrak m,v}^2t^{-1}&\\&1
\end{psmallmatrix}
b'_v
\right)
\chi_v(t)|\!|t|\!|_v^{s-1/2}\,d^\times t,
&v\notin S_0,\\[1em]
\displaystyle
\int_{F_v^\times}
\widetilde W_{\xi_v}
\!\left(
\begin{psmallmatrix}t^{-1}&\\&1\end{psmallmatrix}
b'_vk_v
\right)
\chi_v(t)|\!|t|\!|_v^{s-1/2}\,d^\times t,
&v\in S_0.
\end{cases}
\]
Here $b'_v$ and
\[
k_v=
\begin{psmallmatrix}
(1+{\sf e}_v{\sf e}'_v)/{\sf f}_v&-{\sf e}'_v\\
-{\sf e}_v&{\sf f}_v
\end{psmallmatrix}
\]
are as in \eqref{eqn-feeprime}.

Suppose first that $v\notin S\cup S_0$.
Then
$\tilde\xi_{\omega^{-1},v}=\tilde\xi_v$,
${\sf f}_{\mathfrak m,v}=1$, and ${\sf f}_v=1$.
Since ${\sf e}'_v\in\o_v$, we have
\begin{align*}
{}^tb_v'^{-1}\cdot\tilde\xi_v
&=
\begin{psmallmatrix}1&\\-{\sf e}'_v&1\end{psmallmatrix}
\begin{psmallmatrix}{\sf n}_v&\\&1\end{psmallmatrix}
\cdot\xi_v\\
&=
\begin{psmallmatrix}{\sf n}_v&\\&1\end{psmallmatrix}
\begin{psmallmatrix}1&\\-{\sf n}_v{\sf e}'_v&1\end{psmallmatrix}
\cdot\xi_v
=\tilde\xi_v.
\end{align*}
The last equality follows from the
$K_{1,v}(\p_v^{n(\pi_v)})$-invariance of $\xi_v$.
Consequently,
\begin{align*}
I_v
&=
\int_{F_v^\times}
\widetilde W_{\tilde\xi_v}
\!\left(
\begin{psmallmatrix}{\sf n}_vt^{-1}&\\&1\end{psmallmatrix}
\right)
\chi_v(t)|\!|t|\!|_v^{s-1/2}\,d^\times t\\
&=
\int_{F_v^\times}
\widetilde W_{\xi_v}
\!\left(
\begin{psmallmatrix}t^{-1}&\\&1\end{psmallmatrix}
\right)
\chi_v(t)|\!|t|\!|_v^{s-1/2}\,d^\times t\\
&=
\Psi(1-s;\widetilde W_{\xi_v},\chi_v^{-1}).
\end{align*}
As in \S\ref{CPS}, this vanishes unless $\chi_v$ is
unramified, and in that case it equals
\[
\chi_v({\sf d}_v)^{-1}|{\sf d}_v|_v^{1/2-s}
\epsilon(s,\pi_v\times\chi_v,\psi_v)\,
L(1-s,\widetilde\pi_v\times\chi_v^{-1}).
\]

Next suppose that $v\in S_0$. Then
${\sf n}_v={\sf f}_{\mathfrak m,v}=1$ and
$\tilde\xi_{\omega^{-1},v}=\xi_v$.
The identity
\[
b'_vk_v=wb_vw^{-1}
       ={\sf f}_v\,{}^tb_v^{-1}
\]
gives
\[
\widetilde W_{\xi_v}(g b'_vk_v)
=
\chi_{\pi_v}({\sf f}_v)^{-1}
\widetilde W_{b_v\cdot\xi_v}(g).
\]
Changing variables $y=t^{-1}$ therefore yields
\[
I_v=
\chi_{\pi_v}({\sf f}_v)^{-1}
\Psi\!\left(
1-s;\widetilde W_{b_v\cdot\xi_v},\chi_v^{-1}
\right).
\]
This argument applies to every place in $S_0$ and
requires no compact invariance of $\xi_v$.

Finally, suppose that $v\in S$.
Here ${\sf n}_v={\sf f}_v=1$, and
\[
\tilde\xi_{\omega^{-1},v}
=
\int_{\o_v^\times}
\omega_v(u)^{-1}
\pi_v\!\left(
\begin{psmallmatrix}
1&\\-\varpi_v^{-n_v}u&1
\end{psmallmatrix}
\right)\xi_v\,d^\times u.
\]
The matrix
\[
{}^tb_v'^{-1}
=
\begin{psmallmatrix}1&\\-{\sf e}'_v&1\end{psmallmatrix}
\]
fixes the spherical vector $\xi_v$ and commutes with
the lower unipotent matrices in this integral.
It therefore fixes $\tilde\xi_{\omega^{-1},v}$.
Using the relation between right translation and
Whittaker duality, we obtain
\begin{align*}
I_v
&=
\int_{F_v^\times}\int_{\o_v^\times}
\omega_v(u)^{-1}
\widetilde W_{\xi_v}\!\left(
\begin{psmallmatrix}
{\sf f}_{\mathfrak m,v}^2t^{-1}&\\&1
\end{psmallmatrix}
\begin{psmallmatrix}
1&\varpi_v^{-n_v}u\\&1
\end{psmallmatrix}
\right)
\chi_v(t)|\!|t|\!|_v^{s-1/2}
\,d^\times u\,d^\times t.
\end{align*}

Since $\widetilde W_{\xi_v}$ is right invariant under
$\begin{psmallmatrix}u^{-1}&\\&1\end{psmallmatrix}$
for $u\in\o_v^\times$, the matrix identity used in
\eqref{eqn-double}, followed by the substitution
$t=ru$, gives
\begin{align*}
I_v
&=
\left(
\int_{\o_v^\times}
\omega_v(u)^{-1}\chi_v(u)\,d^\times u
\right)\\
&\qquad\times
\int_{F_v^\times}
\widetilde W_{\xi_v}\!\left(
\begin{psmallmatrix}
{\sf f}_{\mathfrak m,v}^2r^{-1}&\\&1
\end{psmallmatrix}
u_v
\right)
\chi_v(r)|\!|r|\!|_v^{s-1/2}\,d^\times r.
\end{align*}
The unit integral is zero unless
$\chi_v\omega_v^{-1}$ is unramified, in which case
it is $1$.
Under this condition, using
\[
R(u_v)\widetilde W_{\xi_v}
=
\widetilde W_{\bar u_v\cdot\xi_v}
\]
and changing variables
$y={\sf f}_{\mathfrak m,v}^2r^{-1}$, we conclude that
\[
I_v
=
\chi_v(\varpi_v^{2n_v})q_v^{-n_v(2s-1)}
\Psi\!\left(
1-s;\widetilde W_{\bar u_v\cdot\xi_v},\chi_v^{-1}
\right).
\]
Put
\[
C_v(s,\chi_v)
=
\chi_v({\sf d}_v)^{-1}|{\sf d}_v|_v^{1/2-s}.
\]
Arguing as in Lemma~\ref{lemma-m2}, we have 
\begin{align*}
&\Psi\!\left(
1-s;\widetilde W_{\bar u_v\cdot\xi_v},\chi_v^{-1}
\right)\\
&\quad=
C_v(s,\chi_v)
(\chi_{\pi_v}\chi_v^2)({\sf d}_v)
|{\sf d}_v|_v^{2s-1}
\Psi\!\left(
1-s;
\widetilde{L({\sf d}_v)W_{\bar u_v\cdot\xi_v}},
\chi_v^{-1}
\right).
\end{align*}

The local vanishing conditions are together equivalent
to $\chi\sim\omega$. Moreover,
\[
\chi_\pi({\sf f})
=
\prod_{v\in S_0}\chi_{\pi_v}({\sf f}_v),
\qquad
\chi_\pi(\mathfrak m)
=
\prod_{v\in S}\chi_{\pi_v}(\varpi_v^{n_v}).
\]
Thus the factors involving ${\sf f}_v$ cancel, and
multiplying the remaining local factors proves the
stated formula in the domain of absolute convergence.
The identity then holds wherever the Mellin transform
admits meromorphic continuation.
\end{proof}

\begin{remark}\label{rmk-lfe2}
For $v\in S$, the Whittaker functions appearing in
Lemmas~\ref{lemma-m3} and~\ref{lemma-m4} are not,
in general, dual to each other.
Indeed, since left and right translations commute,
the Whittaker function in Lemma~\ref{lemma-m3} is
\[
L({\sf d}_v)W_{u_v\cdot\xi_v}
=
R(u_v)L({\sf d}_v)W_{\xi_v}.
\]
Its dual is
\[
R(\bar u_v)\,
\widetilde{L({\sf d}_v)W_{\xi_v}},
\]
whereas the Whittaker function in
Lemma~\ref{lemma-m4} is
\[
\widetilde{L({\sf d}_v)W_{\bar u_v\cdot\xi_v}}
=
R(u_v)\,
\widetilde{L({\sf d}_v)W_{\xi_v}}.
\]
Thus the local functional equation \eqref{lfe}
does not directly relate the two local integrals.
We derive the required relation below.
\end{remark}

\subsection{The local functional equation for $v\in S$}
Suppose $n_v={\rm ord}_v(\m)$, $u_v=u(\varpi_v^{-n_v})$ and $\bar{u}_v=\bar{u}(-\varpi_v^{-n_v})$ as before. In this section, we only consider characters $\chi_v$ so that $\chi_v\omega_v^{-1}$ is unramified. In particular, this means $n_v\geq m_v={\rm cond}(\chi_v)$. Since $S$ is disjoint from the support of ${\mathfrak n}$, we also have that $\pi_v$ is unramified. 
Put
\[
D_{\m}(s,\pi_v\times\chi_v)
=
\frac{\chi_v(\varpi_v^{m_v-n_v})}
     {\mathscr G(\chi_v,\psi_v,
                 {\sf d}_v^{-1}\varpi_v^{-m_v})}
\frac{\Psi(s;L({\sf d}_v)W_{u_v\cdot\xi_v},\chi_v)}
     {L(s,\pi_v\times\chi_v)}.
\]
Here $L({\sf d}_v)W_{\xi_v}$ is the normalized
essential Whittaker function for the additive
character $x\mapsto\psi_v({\sf d}_v^{-1}x)$,
whose conductor is $\o_v$.

The normalizing factor can also be written as
\[
\frac{\chi_v(\varpi_v^{-n_v})}
{\chi_v(\varpi_v^{-m_v})
 \mathscr G(\chi_v,\psi_v,
           \varpi_v^{-m_v-\delta_v})}.
\]
For $m_v\geq1$, 
\[
\begin{aligned}
&\chi_v(\varpi_v^{-m_v})
\mathscr G(\chi_v,\psi_v,
          \varpi_v^{-m_v-\delta_v})\\
&\qquad=
q_v^{-m_v/2}(1-q_v^{-1})^{-1}
\chi_v({\sf d}_v)
\epsilon\!\left(\tfrac12,\chi_v^{-1},\psi_v\right).
\end{aligned}
\]
When $m_v=0$, the denominator equals $1$.


\begin{lemma}\label{lemma-D}
$D_{\m}(s,\pi_v\times\chi_v)$ belongs to
${\mathbb C}[q_v^{-s}]$.
\end{lemma}

\begin{proof}
We initially work in a right half-plane of
absolute convergence.
\begin{align*}
&\Psi(s;L({\sf d}_v)W_{u_v\cdot\xi_v},\chi_v)\\
&\quad=
\int_{F_v^\times}
(L({\sf d}_v)W_{\xi_v})(a(t))
\psi_v({\sf d}_v^{-1}\varpi_v^{-n_v}t)
\chi_v(t)|t|_v^{s-1/2}\,d^\times t\\
&\quad=
\sum_{r\geq0}
\mathscr G(\chi_v,\psi_v,
          \varpi_v^{r-n_v-\delta_v})
(L({\sf d}_v)W_{\xi_v})(a(\varpi_v^r))
\chi_v(\varpi_v^r)q_v^{-r(s-1/2)}.
\end{align*}
For $r\geq0$, put
\[
\lambda_r
=
q_v^{r/2}
(L({\sf d}_v)W_{\xi_v})(a(\varpi_v^r)),
\]
and set $\lambda_{-1}=0$.
Thus $\lambda_0=1$.

Suppose first that $m_v=0$, and put
$\gamma_v=\chi_v(\varpi_v)$.
Since $x\mapsto\psi_v({\sf d}_v^{-1}x)$ has
conductor $\o_v$, we have
\[
\mathscr G(\chi_v,\psi_v,
          \varpi_v^{r-n_v-\delta_v})
=
\begin{cases}
1,&r\geq n_v,\\
-\dfrac{1}{q_v-1},&r=n_v-1,\\
0,&r<n_v-1.
\end{cases}
\]
Moreover,
\[
L(s,\pi_v\times\chi_v)
=
\sum_{r\geq0}\lambda_r\gamma_v^r q_v^{-rs}.
\]
Consequently,
\begin{align}
&\Psi(s;L({\sf d}_v)W_{u_v\cdot\xi_v},\chi_v)
\notag\\
&\quad=
-\frac{q_v}{q_v-1}
\lambda_{n_v-1}\gamma_v^{n_v-1}
q_v^{-(n_v-1)s}
+
L(s,\pi_v\times\chi_v)
-
\sum_{r=0}^{n_v-2}
\lambda_r\gamma_v^r q_v^{-rs}.
\label{eqn-localint}
\end{align}
The final sum is empty when $n_v=1$.
Since the reciprocal of $L(s,\pi_v\times\chi_v)$
is a polynomial in $q_v^{-s}$, and the
normalizing factor defining $D_{\m}$ is
independent of $s$, this proves the assertion
when $m_v=0$.

Now suppose that $m_v\geq1$.
The Gauss integral
\[
\mathscr G(\chi_v,\psi_v,
          \varpi_v^{r-n_v-\delta_v})
\]
vanishes unless $r=n_v-m_v$.
Hence
\begin{align*}
&\Psi(s;L({\sf d}_v)W_{u_v\cdot\xi_v},\chi_v)\\
&\quad=
\mathscr G(\chi_v,\psi_v,
          \varpi_v^{-m_v-\delta_v})
(L({\sf d}_v)W_{\xi_v})
   (a(\varpi_v^{n_v-m_v}))\\
&\qquad\qquad\times
\chi_v(\varpi_v^{n_v-m_v})
q_v^{-(n_v-m_v)(s-\frac12)}.
\end{align*}
Since $\pi_v$ is unramified and $\chi_v$ is
ramified, $L(s,\pi_v\times\chi_v)=1$.
The normalizing factors therefore cancel, giving
\begin{equation}\label{eqn-Dir2}
D_{\m}(s,\pi_v\times\chi_v)
=
\lambda_{n_v-m_v}q_v^{-(n_v-m_v)s}.
\end{equation}
Since $n_v\geq m_v$, this is a polynomial in
$q_v^{-s}$.
\end{proof}

\begin{lemma}\label{lemma-dirfe}
One has the functional equation 
\[
D_{\m}(s,\pi_v\times\chi_v)=q_{v}^{(n_v-m_v)(1-2s)}\chi_{\pi_v}(\varpi_v^{n_v-m_v})D_{\m}(1-s,\widetilde{\pi}_v\times\chi_v^{-1}).
\]
\end{lemma}
\begin{proof}
Since $\pi_v$ is unramified, the normalized
spherical Whittaker functions satisfy
\[
L({\sf d}_v)W_{\widetilde\pi_v}
=
\widetilde{L({\sf d}_v)W_{\xi_v}}.
\]
Moreover,
\[
(L({\sf d}_v)W_{\widetilde\pi_v})(g)
=
\chi_{\pi_v}(\det g)^{-1}
(L({\sf d}_v)W_{\xi_v})(g).
\]
If $m_v\geq1$, put $k=n_v-m_v$.
By \eqref{eqn-Dir2} and the preceding identity,
\[
D_{\m}(1-s,\widetilde\pi_v\times\chi_v^{-1})
=
\chi_{\pi_v}(\varpi_v^k)^{-1}
\lambda_k q_v^{-k(1-s)}.
\]
Multiplying by
$q_v^{k(1-2s)}\chi_{\pi_v}(\varpi_v^k)$ gives
$\lambda_kq_v^{-ks}
=D_{\m}(s,\pi_v\times\chi_v)$,
as required.

Suppose now that $m_v=0$, and write $n=n_v\geq1$. Let $\{\alpha_v,\beta_v\}$ be the Satake parameters associated with $\pi_v$. With the notation of Lemma~\ref{lemma-D}, we have
\[
\lambda_r=\sum_{j=0}^{r}\alpha_v^{r-j}\beta_v^{j}, r\geq 0.
\] 
We have  
\[
L(s,\pi_v\times\chi_v)=(1-\lambda_1\gamma_vq_v^{-s}+\chi_{\pi_v}(\varpi_v)\gamma_v^2q_v^{-2s})^{-1}\mbox{ with } \chi_{\pi_v}(\varpi_v)=\alpha_v\beta_v,
\]
and the recurrence relation $\lambda_r\lambda_1=\lambda_{r+1}+\chi_{\pi_v}(\varpi_v)\lambda_{r-1}$ for $r\geq 1$. It follows from the latter relation that 
\begin{eqnarray}
\sum_{r=0}^{n-2}\lambda_r\gamma_v^{r}q_v^{-rs}&=&\left(\lambda_{n-2}\gamma_v^n\chi_{\pi_v}(\varpi_v)q_v^{-ns}-\lambda_{n-1}\gamma_v^{n-1}q_v^{-(n-1)s}+1\right)\sum_{r=0}^{\infty}\lambda_r\gamma_v^rq_v^{-rs}\nonumber\\&=&\left(\lambda_{n-2}\gamma_v^n\chi_{\pi_v}(\varpi_v)q_v^{-ns}-\lambda_{n-1}\gamma_v^{n-1}q_v^{-(n-1)s}+1\right)L(s,\pi_v\times\chi_v).\nonumber
\end{eqnarray}
Inserting this into (\ref{eqn-localint}), we obtain
\begin{align}
\Psi(s;L({\sf d}_v)W_{u_v\cdot\xi_v},\chi_v)
&= -\frac{q_v}{q_v-1}\,
\lambda_{n-1}\gamma_v^{n-1}q_v^{-(n-1)s}
\notag\\
&\quad-
L(s,\pi_v\times\chi_v)
\Big(
\lambda_{n-2}\gamma_v^n
\chi_{\pi_v}(\varpi_v)q_v^{-ns}
-
\lambda_{n-1}\gamma_v^{n-1}q_v^{-(n-1)s}
\Big).
\notag
\end{align}

Therefore $D_{\m}(s,\pi_v\times\chi_v)$ equals 
\[
-\frac{q_v}{q_v-1}\lambda_{n-1}\gamma_v^{-1}q_v^{-(n-1)s}\left(1-\lambda_1\gamma_vq_v^{-s}+\chi_{\pi_v}(\varpi_v)\gamma_v^2q_v^{-2s}\right)-\lambda_{n-2}\chi_{\pi_v}(\varpi_v)q_v^{-ns}+\lambda_{n-1}\gamma_v^{-1}q_v^{-(n-1)s}.
\]
Collecting like powers of $q_v^{-s}$ in the above expression and using the recurrence relation we obtain 
\begin{equation}\label{eqn-Dir1}
D_{\m}(s,\pi_v\times\chi_v)=\frac{q_v^{-(n-1)s}}{q_v-1}\left(\lambda_nq_v^{1-s}-\lambda_{n-1}\gamma_v^{-1}+\lambda_{n-2}\chi_{\pi_v}(\varpi_v)q_v^{-s}-\lambda_{n-1}\chi_{\pi_v}(\varpi_v)\gamma_vq_v^{1-2s}\right).
\end{equation}
Put $A_v=\chi_{\pi_v}(\varpi_v)$.
Replacing $\pi_v$ by $\widetilde{\pi}_v$ and $\chi_v$
by $\chi_v^{-1}$ sends
\[
\lambda_r\longmapsto A_v^{-r}\lambda_r,
\qquad
A_v\longmapsto A_v^{-1},
\qquad
\gamma_v\longmapsto\gamma_v^{-1}.
\]
Substituting these into \eqref{eqn-Dir1} gives
\begin{align*}
&D_{\mathfrak m}(1-s,\widetilde{\pi}_v\times\chi_v^{-1})\\
&\quad=
\frac{A_v^{-n}q_v^{(n-1)(s-1)}}{q_v-1}
\left(
\lambda_nq_v^s
-A_v\lambda_{n-1}\gamma_v
+A_v\lambda_{n-2}q_v^{s-1}
-\lambda_{n-1}\gamma_v^{-1}q_v^{2s-1}
\right).
\end{align*}
Multiplying by $A_v^nq_v^{n(1-2s)}$ and comparing with
\eqref{eqn-Dir1} proves the asserted functional equation.
\end{proof}
 
\begin{lemma}\label{lemma-mi2}
 Under the hypothesis of Theorem~\ref{mainthm}, for $\omega\in X(\m)$, $t\in \A_F^{\times}$, and ${\sf f},{\sf e},{\sf e}'$ as in Lemma~\ref{lemma-mi1}, we have 
 \begin{equation}\label{eqn-twist}
 \phi_{\omega}\left(\begin{psmallmatrix}t{\sf f}&t{\sf e}\\&1\end{psmallmatrix}\right)=\epsilon_{[\omega]}\omega_S(-1)\chi_{\pi}(\m)\chi_{\pi}({\sf f})(R(k_{S_0})\widetilde{\phi}_{\omega^{-1}})\left(\begin{psmallmatrix}{\sf n}{\sf f}_{\m}^{2}t^{-1}{\sf f}&{\sf n}{\sf f}_{\m}^{2}t^{-1}{\sf e}'\\&1\end{psmallmatrix}\right).
 \end{equation}
 \end{lemma}
 \begin{proof}
Let $\chi$ be a unitary id\`ele class character.
By Lemmas~\ref{lemma-m3} and~\ref{lemma-m4}, both Mellin
transforms vanish unless $\chi\sim\omega$.
Suppose that $\chi\sim\omega$.

For $v\in S$, put $m_v={\rm cond}(\chi_v)$.
The definition of $D_{\mathfrak m}$ and
Lemma~\ref{lemma-dirfe} give
\begin{align*}
&(\chi_{\pi_v}\chi_v^2)(\varpi_v^{n_v})
q_v^{-n_v(2s-1)}
(\chi_{\pi_v}\chi_v^2)({\sf d}_v)
|{\sf d}_v|_v^{2s-1}\\
&\qquad\times
\frac{
\Psi\!\left(
1-s;
\widetilde{L({\sf d}_v)W_{\bar u_v\cdot\xi_v}},
\chi_v^{-1}
\right)}
{L(1-s,\widetilde\pi_v\times\chi_v^{-1})}\\
&\quad=
\Xi_v(s,\pi_v\times\chi_v,\psi_v)
\frac{
\Psi\!\left(
s;L({\sf d}_v)W_{u_v\cdot\xi_v},\chi_v
\right)}
{L(s,\pi_v\times\chi_v)},
\end{align*}
where
\[
\begin{aligned}
\Xi_v(s,\pi_v\times\chi_v,\psi_v)
&=
(\chi_{\pi_v}\chi_v^2)({\sf d}_v)
|{\sf d}_v|_v^{2s-1}\\
&\quad\times
(\chi_{\pi_v}\chi_v^2)(\varpi_v^{m_v})
q_v^{-m_v(2s-1)}
\frac{
\mathscr G(\chi_v^{-1},\psi_v,
          \varpi_v^{-m_v-\delta_v})
}{
\mathscr G(\chi_v,\psi_v,
          \varpi_v^{-m_v-\delta_v})
}.
\end{aligned}
\]
Here we have used
\[
\widetilde{L({\sf d}_v)W_{\xi_v}}
=
L({\sf d}_v)W_{\widetilde\pi_v}
\]
and
\[
\widetilde{L({\sf d}_v)W_{\bar u_v\cdot\xi_v}}
=
R(u_v)\widetilde{L({\sf d}_v)W_{\xi_v}}.
\]

If $m_v=0$, both Gauss integrals are $1$.
Since $\pi_v$ and $\chi_v$ are unramified,
\[
\epsilon(s,\pi_v\times\chi_v,\psi_v)
=
(\chi_{\pi_v}\chi_v^2)({\sf d}_v)
|{\sf d}_v|_v^{2s-1}.
\]
Also $\chi_v(-1)=1$, so
\[
\Xi_v(s,\pi_v\times\chi_v,\psi_v)
=
\epsilon(s,\pi_v\times\chi_v,\psi_v)\chi_v(-1).
\]

If $m_v\geq1$, formulas~\eqref{eqn-gaussrel}
and~\eqref{eqn-tate}, together with the multiplicativity
of local $\epsilon$-factors for the unramified
representation $\pi_v$ give
\[
\begin{aligned}
&(\chi_{\pi_v}\chi_v^2)(\varpi_v^{m_v})
q_v^{-m_v(2s-1)}
\frac{
\mathscr G(\chi_v^{-1},\psi_v,
          {\sf d}_v^{-1}\varpi_v^{-m_v})
}{
\mathscr G(\chi_v,\psi_v,
          {\sf d}_v^{-1}\varpi_v^{-m_v})
}\\
&\qquad=
\epsilon\!\left(
s,\pi_v\times\chi_v,
\psi_v^{\circ}
\right)\chi_v(-1).
\end{aligned}
\]
Since
\[
\begin{aligned}
\epsilon(s,\pi_v\times\chi_v,\psi_v)
&=
(\chi_{\pi_v}\chi_v^2)({\sf d}_v)
|{\sf d}_v|_v^{2s-1}\\
&\qquad\times
\epsilon\!\left(
s,\pi_v\times\chi_v,
\psi_v({\sf d}_v^{-1}\,\cdot)
\right),
\end{aligned}
\]
it follows that 
\[
\Xi_v(s,\pi_v\times\chi_v,\psi_v)
=
\epsilon(s,\pi_v\times\chi_v,\psi_v)\chi_v(-1)
\]
in both cases.

At each $v\in S_0$, the local functional equation
\eqref{lfe} gives
\[
\frac{
\Psi(1-s;\widetilde W_{b_v\cdot\xi_v},\chi_v^{-1})
}{
L(1-s,\widetilde\pi_v\times\chi_v^{-1})
}
=
\epsilon(s,\pi_v\times\chi_v,\psi_v)
\frac{
\Psi(s;W_{b_v\cdot\xi_v},\chi_v)
}{
L(s,\pi_v\times\chi_v)
}.
\]
Combining these identities with
Lemmas~\ref{lemma-m3} and~\ref{lemma-m4},
and noting that the same factor
$C_\psi^{S_0}(s,\chi)$ occurs in both Mellin
transforms, we obtain
\begin{align*}
\frac{
{\mathscr M}(s,\widetilde f_{[\omega]},\chi)
}{
\Lambda(1-s,\widetilde\pi\times\chi^{-1})
}
&=
\left(\prod_{v\in S}\chi_v(-1)\right)
\epsilon(s,\pi\times\chi,\psi)
\frac{
{\mathscr M}(s,f_{[\omega]},\chi)
}{
\Lambda(s,\pi\times\chi)
}.
\end{align*}
These are identities of meromorphic functions.

Since $\chi\sim\omega$, we have
\[
\epsilon_\chi=\epsilon_{[\omega]},
\qquad
\prod_{v\in S}\chi_v(-1)=\omega_S(-1).
\]
Using the assumed global functional equation
\eqref{funceq} and $\omega_S(-1)^2=1$, it follows that
\[
{\mathscr M}(s,f_{[\omega]},\chi)
=
\epsilon_{[\omega]}\omega_S(-1)
{\mathscr M}(s,\widetilde f_{[\omega]},\chi).
\]
Together with the vanishing when $\chi\not\sim\omega$,
this gives the identity for every unitary id\`ele
class character $\chi$.
Applying Lemma~\ref{lemma-JL} yields \eqref{eqn-twist}.
\end{proof}
 
We now make the normalization described in
Remark~(\ref{rmk-normalization}).
Put $c=\epsilon_{[1]}\ne0$.
We replace each auxiliary dual function by $c$ times
that function and each $\epsilon_\chi$
by $\epsilon_\chi/c$.
We retain the same notation for the normalized
functions and multipliers; thus $\epsilon_{[1]}=1$
from this point onward.
The local Whittaker functions and local
$\epsilon$-factors remain unchanged.

\section{Proof of Theorem~\ref{mainthm}}
We retain ${\sf f},{\sf e},{\sf e}'$ as in the previous lemma, and put
\[
g=\begin{pmatrix} t{\sf f} & t{\sf e} \\ 0 & 1 \end{pmatrix}.
\]
For an admissible $\m$ as before, let $S={\rm Supp}(\m)$. Then $S\cap S_0=\emptyset$. Let 
\[
G=X(\m),
\qquad
H=X(1).
\]
With $n_v=\operatorname{ord}_v(\m)$, let
\[
 R_{\m}=\prod_{v\in S}\o_v/\p_v^{n_v},
 \quad {\mathcal K}_{\m}=R_{\m}^{\times},
 \quad N_{\m}=|{\mathcal K}_{\m}|,\qquad
\Psi_{\m}(x)
=
\prod_{v\in S}
\psi_v(\varpi_v^{-n_v-\delta_v}x_v).
\]
The group ${\mathcal K}_{\m}$ may also be realized as 
\[
{\mathcal K}_{\m}=\prod_{v\in S}{\mathcal K}_v, 
\]
where ${\mathcal K}_v=\o_v^{\times}/(1+\p_v^{n_v})$. 

Let $\widehat{\mathcal K}_{\m}=\prod_{v\mid \m}\widehat{\mathcal K}_v$ denote the corresponding character group. Restriction to the local unit group defines a surjective homomorphism from 
\[
G\longrightarrow \widehat{\mathcal K}_{\m},
\qquad
\omega\longmapsto\omega_S,
\]
where
\[
\omega_S((u_v)_{v\in S})
=
\prod_{v\in S}\omega_v(u_v).
\]
Its kernel is precisely $H$ and consequently 
\[
G/H \simeq \widehat{\mathcal K}_{\m}.
\]

For each $\omega_S=\prod_{v\in S}\omega_v\in \widehat{\mathcal K}_{\m}$, choose a representative
$\omega\in G$. The functions $\phi_{[\omega]}$ and
$\widetilde\phi_{[\omega]}$, as well as the multiplier
$\epsilon_{[\omega]}$, are independent of this choice.
Take any finite id\`{e}le $b=(b_v)$ so that $b_v\in \o_v^{\times}$ for all $v\in S$, and $b_v=1$ otherwise, multiply both sides of (\ref{eqn-twist}) by $\omega(b)$ and sum over all $\omega\in G/H$. One checks that
\begin{eqnarray}
\sum_{\omega\in G/H}\omega(b)\phi_{[\omega]}(g)&=&\sum_{\omega\in G/H}\sum_{\gamma\in F^{\times}}\int\limits_{\prod_{v\in S}\o_v^{\times}} W_{\xi}\left(\begin{psmallmatrix}\gamma &\\&1\end{psmallmatrix}g\begin{psmallmatrix}1&{\sf w}u\\&1\end{psmallmatrix}\right)\omega(ub)d^{\times}u\nonumber\\&=&\sum_{\gamma\in F^{\times}}
\prod_{v'\notin S}W_{\xi_{v'}
}\left(\begin{psmallmatrix}\gamma&\\&1\end{psmallmatrix}g_{v'}\right)\times \prod_{v\in S}I_v(\gamma,g_v),\label{eqn-Iv}
\end{eqnarray}
where 
\[
I_v(\gamma,g_v)=\int\limits_{\o_v^{\times}} W_{\xi_v}\left(\begin{psmallmatrix}\gamma &\\&1\end{psmallmatrix}g_v\begin{psmallmatrix}1&{\sf w}_vub_v^{-1}\\&1\end{psmallmatrix}\right)(\sum_{\omega_v\in \widehat{\mathcal K}_v}\omega_v(u))d^{\times}u.
\]
But $\sum_{\omega_v\in \widehat{\mathcal K}_v}\omega_v(u_{v})=0$ unless $u_v\in U_v^{(n_v)}$. Consequently, for $v\in S$, 
\[
I_v(\gamma,g_v)= W_{\xi_v}\left(\begin{psmallmatrix}\gamma &\\&1\end{psmallmatrix}g_v\begin{psmallmatrix}1&{\sf w}_vb_v^{-1}\\&1\end{psmallmatrix}\right).
\]
Inserting this into equation (\ref{eqn-Iv}) and using the definition of $\phi$, we obtain 
\begin{equation}\label{eqn-LHS}
\sum_{\omega}\omega(b)\phi_{[\omega]}\left(\begin{psmallmatrix}t{\sf f}&t{\sf e}\\&1\end{psmallmatrix}\right)=\phi\left(\begin{psmallmatrix}t{\sf f}&t{\sf e}\\&1\end{psmallmatrix}\begin{psmallmatrix}1&{\sf w}b^{-1}\\&1\end{psmallmatrix}\right).
\end{equation}

On the other hand, the right hand side of (\ref{eqn-twist}) after summation is given by 
\begin{equation}
\chi_{\pi}({\sf f})\chi_{\pi}(\m)
\!\!\int\limits_{\prod_{v\in S}\o_v^{\times}}\!
\Biggl(
\sum_{\omega\in G/H}\epsilon_{[\omega]}\,\omega^{-1}(u)
\Biggr)
(R(k_{S_0})\widetilde{\phi})\!\left(
\begin{psmallmatrix}
{\sf n}{\sf f}_{\m}^{2}t^{-1}{\sf f} & {\sf n}{\sf f}_{\m}^{2}t^{-1}{\sf e}'\\
& 1
\end{psmallmatrix}
\begin{psmallmatrix}
1 & -{\sf w}ub\\
& 1
\end{psmallmatrix}
\right)
\,d^{\times}u.\nonumber
\end{equation}
Inserting the Fourier expansion for $\widetilde{\phi}$ (cf. (\ref{dual-FE})) into the above equation results in
\begin{align}\label{eqn-RHS}
&\chi_{\pi}({\sf f})\,\chi_{\pi}(\m)
\sum_{\gamma\in F^{\times}}
\widetilde{W}_{S_{0}}\!\left(
\begin{psmallmatrix}
\gamma t_{S_0}^{-1}&\\&1
\end{psmallmatrix}
b_{S_0}'k_{S_0}
\right)
\tilde{c}^{S_0}(\gamma {\sf n}{\sf f}^{2}_{\m}{\sf f}t^{-1})
\\[3pt]
&\quad\times
\left(
\int\limits_{\prod_{v\in S}\o_v^{\times}}
\Biggl(
\sum_{\omega_S\in \widehat{\mathcal K}_{\m}}
\epsilon_{[\omega]}\,\omega_S^{-1}(u)
\Biggr)
\Psi_{\m}\!\bigl(
-{\sf d}\gamma{\sf n}{\sf f}_{\m}^{2}
{\sf f}t^{-1}(u-1)b
\bigr).
d^{\times}u
\right)
\nonumber\\
&\quad\times
\psi\!\bigl(-\gamma{\sf n}{\sf f}^{2}_{\m}{\sf f}t^{-1}{\sf w}b\bigr)
\psi^{S_0}\!\bigl(\gamma {\sf n}{\sf f}^{2}_{\m}t^{-1}{\sf e}'\bigr).
\nonumber
\end{align}
For $x\in R_{\m}$, 
put 
\[
S_{\m}(x)=\int\limits_{\prod_{v\in S}\o_v^{\times}}\hat{\epsilon}_{\m}(u)\Psi_{\m}(-(u-1)x)d^{\times}u,
\]
where $\hat{\epsilon}_{\m}(u)=\sum_{\omega_S\in {\widehat K}_{\m}}\epsilon_{[\omega]}\omega_S^{-1}(u)$. It becomes a finite sum 
\[
S_{\m}(x)=\frac{1}{N_{\m}}\sum_{u\in {\mathcal K}_{\m}}\hat{\epsilon}_{\m}(u)\Psi_{\m}((1-u)x).
\]
Now, we can express (\ref{eqn-RHS}) as 
\begin{align}\label{eqn-RHS2}
&\chi_{\pi}({\sf f})\,\chi_{\pi}(\m)
\sum_{\gamma\in F^{\times}}
\widetilde{W}_{S_0}\!\left(
\begin{psmallmatrix}
\gamma t_{S_0}^{-1} & \\
& 1
\end{psmallmatrix}
b_{S_0}'k_{S_0}
\right)
\tilde{c}^{S_0}(\gamma {\sf n}{\sf f}^{2}_{\m}{\sf f}t^{-1})
\\
&\quad\times
S_{\m}\!\bigl(
{\sf d}\gamma{\sf n}{\sf f}_{\m}^{2}
{\sf f}t^{-1}b
\bigr)\,
\psi\!\bigl(-\gamma{\sf n}{\sf f}^{2}_{\m}{\sf f}t^{-1}{\sf w}b\bigr)\,
\psi^{S_0}\!\bigl(\gamma {\sf n}{\sf f}^{2}_{\m}t^{-1}{\sf e}'\bigr).
\nonumber
\end{align}

We absorb the scalar $c$ used to normalize the
auxiliary dual function into its coefficients
$\tilde c^{S_0}$, retaining the same notation.
The local Whittaker functions remain unchanged. Now,
equating (\ref{eqn-LHS}) and (\ref{eqn-RHS2}), we obtain
\begin{align}\label{eqn-R}
\phi\!\left(
\begin{psmallmatrix}
t {\sf f} & t {\sf e}\\
& 1
\end{psmallmatrix}
\begin{psmallmatrix}
1 & {\sf w}b^{-1}\\
& 1
\end{psmallmatrix}
\right)
&=
\chi_{\pi}({\sf f})\,\chi_{\pi}(\m)
\sum_{\gamma\in F^{\times}}
\widetilde{W}_{S_0}\!\left(
\begin{psmallmatrix}
\gamma t_{S_0}^{-1} & \\
& 1
\end{psmallmatrix}
b_{S_0}'k_{S_0}
\right)
\tilde{c}^{S_0}(\gamma {\sf n}{\sf f}^{2}_{\m}{\sf f}t^{-1})
\\
&\quad\times
S_{\m}\!\bigl(
{\sf d}\gamma{\sf n}{\sf f}_{\m}^{2}
{\sf f}t^{-1}b
\bigr)\,
\psi\!\bigl(-\gamma{\sf n}{\sf f}^{2}_{\m}{\sf f}t^{-1}{\sf w}b\bigr)\,
\psi^{S_0}\!\bigl(
\gamma{\sf n}{\sf f}_{\m}^{2}t^{-1}{\sf e}'
\bigr).
\nonumber
\end{align}
Only those $\gamma$ for which
\[
\operatorname{div}^{S_0}
\bigl(
{\sf d}\gamma{\sf n}{\sf f}_{\m}^{2}
{\sf f}t^{-1}
\bigr)\geq0
\]
can contribute, by the support of
$\tilde c^{S_0}$.
For these terms, the argument
\[
{\sf d}\gamma{\sf n}{\sf f}_{\m}^{2}
{\sf f}t^{-1}b
\]
defines an element of $R_{\m}$.

\subsection{A special case of the automorphy relation}\label{ss-aut}
Fix an effective divisor $\mathfrak H$, supported
outside $S_0$ and coprime to $\mathfrak n$.
Choose the id\`ele ${\sf h}$ with components
\[
{\sf h}_v=
\begin{cases}
\varpi_v^{\operatorname{ord}_v(\mathfrak H)},
    &v\notin S_0,\\
1,&v\in S_0.
\end{cases}
\]
Let $\mathfrak m$ be an admissible modulus
coprime to $\mathfrak H$, and put
$S=\operatorname{Supp}(\mathfrak m)$.

At every nonarchimedean place $v$, set
\[
{\sf f}_{1,v}={\sf f}_{\mathfrak m,v},
\qquad
{\sf e}_{1,v}={\sf h}_v,
\qquad
t_{1,v}={\sf h}_v^{-1},
\]
and
\[
{\sf e}'_{1,v}=
\begin{cases}
-1,&v\in S,\\
0,&v\notin S.
\end{cases}
\]
At each archimedean place $w$, set
\[
{\sf f}_{1,w}={\sf e}_{1,w}=\frac1{\sqrt2},
\qquad
{\sf e}'_{1,w}=-\frac1{\sqrt2},
\qquad
t_{1,w}=\sqrt2.
\]
Then $({\sf f}_1,{\sf e}_1)$ is a reduced pair
of conductor $\mathfrak m$, and ${\sf e}'_1$
satisfies \eqref{eqn-edisjoint}.
Moreover,
\[
\begin{psmallmatrix}
t_1{\sf f}_1&t_1{\sf e}_1\\&1
\end{psmallmatrix}
=
\begin{psmallmatrix}
{\sf f}_{\mathfrak m}{\sf h}^{-1}&1\\&1
\end{psmallmatrix}.
\]



\begin{lemma}\label{lemma-density}
Retain $\mathfrak H$, ${\sf h}$ and $\mathfrak m$
as above. Suppose that there exists $M\in F^\times$
such that
\[
\operatorname{div}^{S_0}(M)=\mathfrak m,
\qquad
M\equiv1\pmod{\mathfrak n\mathfrak H}.
\]
Put
\[
b'=u(-1)a({\sf n}{\sf h}{\sf f}_{\mathfrak m}^{-1}).
\]
Then, for every choice of smooth vectors $\xi_{S_0}$,
\[
\widetilde\phi_{\xi_{S_0}}(wb')
=
\widetilde\phi_{\xi_{S_0}}(b').
\]
\end{lemma}

\begin{proof}
The assertion is unchanged by multiplying the dual function
by a nonzero scalar. We therefore work with the original,
unscaled dual function
$\widetilde{\phi}=\widetilde\varphi_\xi$, and write $\epsilon$
for the original multiplier $\epsilon_{[1]}$.
By Lemma~\ref{lemma-mi1}, we have
\[
\phi(a({\sf n}t^{-1}))
=
\epsilon\widetilde{\phi}(a(t)),
\]
which in turn, using
\[
wa(t)
\begin{psmallmatrix}0&1\\-{\sf n}&0\end{psmallmatrix}
=
-t\,a({\sf n}t^{-1}),
\]
gives
\[
\widetilde{\phi}(wa(t))
=
\varphi_\xi(a({\sf n}t^{-1}))
=
\epsilon\widetilde{\phi}(a(t)),
\qquad t\in\A_F^\times.
\]
For $h\in G_{S_0}$, right translation of
$\widetilde\varphi_\xi$ by $h$ corresponds to
replacing $\xi_{S_0}$ by
$\pi_{S_0}({}^th^{-1})\xi_{S_0}$.
Here ${\sf n}_{S_0}=1$, and the vectors outside
$S_0$ remain unchanged. It follows that
\begin{equation}\label{eqn-density-diagonal}
\widetilde{\phi}(wa(t)h)
=
\epsilon\widetilde{\phi}(a(t)h)
\qquad
t\in\A_F^\times,\quad h\in G_{S_0}.
\end{equation}

Write $A=\o_F^{S_0}$ and put
\[
\mathcal F(h)
=
\widetilde{\phi}(a({\sf n}{\sf h})h),
\qquad h\in G_{S_0}.
\]
We claim that
\begin{align}
\mathcal F(u(r)_{S_0}h)
&=\mathcal F(h),
&&r\in\mathfrak n\mathfrak H A,
\label{eqn-density-upper}\\
\mathcal F(\bar u(s)_{S_0}h)
&=\mathcal F(h),
&&s\in A.
\label{eqn-density-lower}
\end{align}

Since $({\sf n}{\sf h})_{S_0}=1$, left invariance
under $u(r)\in B(F)$ and right invariance under
$K_1^{S_0}(\mathfrak n)$ give
\[
\begin{aligned}
\widetilde{\phi}(a({\sf n}{\sf h})u(r)_{S_0}h)
&=
\widetilde{\phi}(u(r)_{S_0}a({\sf n}{\sf h})h)\\
&=
\widetilde{\phi}(u(-r)^{S_0}a({\sf n}{\sf h})h)\\
&=
\widetilde{\phi}
\left(
a({\sf n}{\sf h})
u(-r/({\sf n}{\sf h}))^{S_0}h
\right)\\
&=
\widetilde{\phi}(a({\sf n}{\sf h})h),
\end{aligned}
\]
proving \eqref{eqn-density-upper}. Here
\[
u(-r/({\sf n}{\sf h}))^{S_0}
\in K_1^{S_0}(\mathfrak n).
\]

Similarly, using \eqref{eqn-density-diagonal},
\[
\begin{aligned}
\epsilon\mathcal F(\bar u(s)_{S_0}h)
&=
\widetilde{\phi}
(wa({\sf n}{\sf h})\bar u(s)_{S_0}h)\\
&=
\widetilde{\phi}
(u(-s)_{S_0}wa({\sf n}{\sf h})h)\\
&=
\widetilde{\phi}
(u(s)^{S_0}wa({\sf n}{\sf h})h)\\
&=
\widetilde{\phi}
(wa({\sf n}{\sf h})\bar u(-s{\sf n}{\sf h})^{S_0}h)\\
&=
\widetilde{\phi}(wa({\sf n}{\sf h})h)
=
\epsilon\mathcal F(h).
\end{aligned}
\]
Since $\mathfrak H$ is effective,
\[
\bar u(-s{\sf n}{\sf h})^{S_0}
\in K_1^{S_0}(\mathfrak n).
\]
Together with $\epsilon\ne0$, this proves
\eqref{eqn-density-lower}.

Now put $h=a(M_{S_0})$ and define the id\`ele
$\delta$ by
\[
\delta_v=
\begin{cases}
M_v{\sf f}_{\mathfrak m,v}^{-1},&v\notin S_0,\\
1,&v\in S_0.
\end{cases}
\]
By hypothesis, $\delta_v\in\o_v^\times$ for
$v\notin S_0$, and
\[
a({\sf n}{\sf h}{\sf f}_{\mathfrak m}^{-1})
=
a({\sf n}{\sf h}/M)h\,a(\delta).
\]
Consequently,
\[
\begin{aligned}
\widetilde\phi(b')
&=
\widetilde\phi
\left(a({\sf n}{\sf h}{\sf f}_{\mathfrak m}^{-1})\right)\\
&=
\widetilde\phi
\left(a({\sf n}{\sf h}/M)h\,a(\delta)\right)\\
&=
\widetilde\phi(a({\sf n}{\sf h})h)
=
\mathcal F(h).
\end{aligned}
\]
On the other hand, the identity
\[
wu(-1)a(M^{-1})
=
M^{-1}a(M)wu(-M)
\]
gives
\[
\widetilde{\phi}(wb')
=
\widetilde{\phi}(wu(-M)a({\sf n}{\sf h})h).
\]
Put
\[
k_M=
\begin{psmallmatrix}
1&(1-M)/({\sf n}{\sf h})\\
-{\sf n}{\sf h}&M
\end{psmallmatrix}.
\]
Direct multiplication shows that
\[
wu(-M)a({\sf n}{\sf h})
=
u(1)a({\sf n}{\sf h})k_M.
\]
Our assumptions on $M$ imply
$k_M^{S_0}\in K_1^{S_0}(\mathfrak n)$, while
\[
(k_M)_{S_0}
=
\bar u(-1)_{S_0}u(1-M)_{S_0}.
\]
As $-1\in A$ and $1-M\in\mathfrak n\mathfrak H A$,
\eqref{eqn-density-upper} and
\eqref{eqn-density-lower} now give
\[
\widetilde{\phi}(wb')
=
\mathcal F((k_M)_{S_0}h)
=
\mathcal F(h)
=
\widetilde{\phi}(b').
\]
\end{proof}

\begin{proposition}\label{prop-special}
Let $\pi$ be as in Theorem~\ref{mainthm}, and retain
$\mathfrak H$ and ${\sf h}$ as above.
Let $\mathfrak m$ be an admissible modulus coprime
to $\mathfrak H$ for which there exists $M\in F^\times$
satisfying
\[
\operatorname{div}^{S_0}(M)=\mathfrak m,
\qquad
M\equiv1\pmod{\mathfrak n\mathfrak H}.
\]
With ${\sf f}_1,{\sf e}_1,{\sf e}'_1$ and $t_1$
chosen as above, we have
\[
\phi_{\xi_{S_0}}\left(
\begin{psmallmatrix}
t_1{\sf f}_1&t_1{\sf e}_1\\
&1
\end{psmallmatrix}
\right)
=
\chi_\pi({\sf f}_1)
(R(k_{S_0})\widetilde\phi_{\xi_{S_0}})\left(
\begin{psmallmatrix}
{\sf n}t_1^{-1}{\sf f}_1&
{\sf n}t_1^{-1}{\sf e}'_1\\
&1
\end{psmallmatrix}
\right),
\]
where $k_{S_0}=k_{S_0}({\sf f}_1,{\sf e}_1)$,
for every choice of smooth vectors $\xi_{S_0}$.
\end{proposition}

\begin{proof}
Write
$(\phi,\widetilde\phi)
=(\phi_{\xi_{S_0}},\widetilde\phi_{\xi_{S_0}})$,
with the normalization fixed above. By construction,
\[
\begin{psmallmatrix}
t_1{\sf f}_1&t_1{\sf e}_1\\
&1
\end{psmallmatrix}
=
\begin{psmallmatrix}
{\sf f}_{\mathfrak m}{\sf h}^{-1}&1\\
&1
\end{psmallmatrix}.
\]
Left ${\rm B}_1(F)$-invariance and
Lemma~\ref{lemma-mi1}, applied to the diagonal matrix,
therefore give
\begin{align*}
\phi\left(
\begin{psmallmatrix}
t_1{\sf f}_1&t_1{\sf e}_1\\
&1
\end{psmallmatrix}
\right)
&=
\phi\left(
\begin{psmallmatrix}
{\sf f}_{\mathfrak m}{\sf h}^{-1}&0\\
&1
\end{psmallmatrix}
\right)\\
&=
\widetilde\phi\left(
\begin{psmallmatrix}
{\sf n}{\sf h}{\sf f}_{\mathfrak m}^{-1}&0\\
&1
\end{psmallmatrix}
\right).
\end{align*}

To evaluate the right-hand side of the asserted identity,
put
\[
B=
\begin{psmallmatrix}
{\sf n}t_1^{-1}{\sf f}_1&
{\sf n}t_1^{-1}{\sf e}'_1\\
&1
\end{psmallmatrix},
\qquad
{\sf f}={\sf n}t_1^{-1},
\qquad
{\sf e}={\sf f}_1^{-1}
       ({\sf n}t_1^{-1}{\sf e}'_1+1).
\]
Since $R(k_{S_0})\widetilde\phi$ is left
${\rm B}_1(F)$-invariant, we have
\begin{equation}\label{eqn-rhs}
(R(k_{S_0})\widetilde\phi)(B)
=
(R(k_{S_0})\widetilde\phi)(b),
\qquad
b=
\begin{psmallmatrix}
{\sf f}_1{\sf f}&{\sf f}_1{\sf e}\\
&1
\end{psmallmatrix}
=u(1)B.
\end{equation}
Our choices give
\[
{\sf f}_w={\sf e}_w=\frac1{\sqrt2}
\quad(w\mid\infty),
\]
and, at finite places,
\[
{\sf f}_v={\sf n}_v{\sf h}_v,
\qquad
{\sf e}_v=
\begin{cases}
0,&v\in S,\\
1,&v\notin S.
\end{cases}
\]
Here ${\sf n}_v={\sf h}_v=1$ for $v\in S$.

Choose ${\sf e}'$ by
\[
{\sf e}'_v=
\begin{cases}
-1/\sqrt2,&v\mid\infty,\\
-{\sf f}_{\mathfrak m,v},&v\in S,\\
-1,&v<\infty,\ v\notin S.
\end{cases}
\]
Define
\[
b'=
\begin{psmallmatrix}
{\sf f}_1^{-1}{\sf f}&
{\sf f}_1^{-1}{\sf e}'\\
&1
\end{psmallmatrix},
\qquad
z=
\begin{psmallmatrix}
{\sf f}_1&0\\
0&{\sf f}_1
\end{psmallmatrix},
\qquad
k'=
\begin{psmallmatrix}
-{\sf e}'&-(1+{\sf e}{\sf e}')/{\sf f}\\
{\sf f}&{\sf e}
\end{psmallmatrix}.
\]
Then
\begin{equation}\label{eqn-decomp}
b=wb'k'z.
\end{equation}
Moreover, the local choices above give
\[
k'^{S_0}\in K_1^{S_0}(\mathfrak n),
\qquad
k'_{S_0}=k_{S_0}^{-1},
\]
and
\[
b'
=
u(-1)a({\sf n}{\sf h}{\sf f}_{\mathfrak m}^{-1}).
\]
Indeed, at $v\notin S\cup S_0$,
\[
k'_v=
\begin{psmallmatrix}
1&0\\
{\sf n}_v{\sf h}_v&1
\end{psmallmatrix}
\in K_{1,v}(\p_v^{n(\pi_v)}),
\]
while at $v\in S$,
\[
k'_v=
\begin{psmallmatrix}
{\sf f}_{\mathfrak m,v}&-1\\
1&0
\end{psmallmatrix}
\in\GL_2(\o_v).
\]

Using \eqref{eqn-decomp}, we obtain
\begin{align*}
(R(k_{S_0})\widetilde\phi)(b)
&=\widetilde\phi(wb'k'zk_{S_0})\\
&=\chi_\pi({\sf f}_1)^{-1}
  \widetilde\phi(wb'k'^{S_0})\\
&=\chi_\pi({\sf f}_1)^{-1}\widetilde\phi(wb')\\
&=\chi_\pi({\sf f}_1)^{-1}\widetilde\phi(b'),
\end{align*}
where the last equality follows from
Lemma~\ref{lemma-density}.

Finally, left invariance under $u(1)\in{\rm B}_1(F)$
gives
\[
\widetilde\phi(b')
=
\widetilde\phi(u(1)b')
=
\widetilde\phi\left(
\begin{psmallmatrix}
{\sf n}{\sf h}{\sf f}_{\mathfrak m}^{-1}&\\
&1
\end{psmallmatrix}
\right).
\]
Combining this with \eqref{eqn-rhs} and the first
displayed identity proves the proposition.
\end{proof}

\begin{lemma}\label{lemma-isolation}
For any irreducible admissible generic
representation $\pi$ of $\GL_2(\A_F)$, let $\mathcal L$
be a nonzero fractional ideal of $A=\o_F^{S_0}$.
Put $\Gamma=\mathcal L\setminus\{0\}$.
For $\gamma_0\in\Gamma$ and $g_{S_0}\in G_{S_0}$, there exists a pure tensor $\xi_{S_0}$
of smooth vectors such that
\[
W_{\xi_{S_0}}(a(\gamma)g_{S_0})
=
\begin{cases}
1,&\gamma=\gamma_0,\\
0,&\gamma\in\Gamma\setminus\{\gamma_0\}.
\end{cases}
\]
\end{lemma}
\begin{proof}
By the product
formula, the diagonal image of $\mathcal L$ in $F_{S_0}=\prod_{v\in S_0}F_v$ is discrete. Fix $\gamma_0\in\Gamma$. We may choose relatively
compact open neighborhoods $V_v$ of $\gamma_0$
in $F_v$, with $0\notin V_v$, such that
\[
\left(\prod_{v\in S_0}V_v\right)
\cap\Gamma=\{\gamma_0\}.
\]
Choose Schwartz--Bruhat functions $h_v$ on $F_v$,
supported in $V_v$, with $h_v(\gamma_0)=1$, and
put $h=\bigotimes_{v\in S_0}h_v$. Thus
\[
h(\gamma)=
\begin{cases}
1,&\gamma=\gamma_0,\\
0,&\gamma\in\Gamma\setminus\{\gamma_0\}.
\end{cases}
\]

Let
\[
\kappa(x)
=
\int_{F_{S_0}}h(y)\psi_{S_0}(-xy)\,dy,
\]
be its Fourier transform with respect to the self-dual additive Haar measure $dy$.
Then $\kappa$ is also a Schwartz--Bruhat function,
and Fourier inversion gives
\[
h(y)=
\int_{F_{S_0}}\kappa(x)\psi_{S_0}(xy)\,dx.
\]
Choose
\[
W'_v\in\mathcal W(\pi_v,\psi_v),
\qquad W'_v(a(\gamma_0))=1
\quad(v\in S_0),
\]
and put $W'=\bigotimes_{v\in S_0}W'_v$.
Let 
\[
W'_h(g)=
\int_{F_{S_0}}\kappa(x)W'(g\,u(x))\,dx.
\]
This defines an element of
$\mathcal W(\pi_{S_0},\psi_{S_0})$.
Indeed, at nonarchimedean places it is a finite
linear combination of right translates of $W'$. At archimedean places, moderate growth of the
smooth globalization and rapid decay of $\kappa$
imply convergence in the smooth Fr\'echet topology \cite[\S 2.5]{Be-K}. Using $a(y)u(x)=u(yx)a(y)$, we obtain
\[
\begin{aligned}
W'_h(a(y))
&=
\int_{F_{S_0}}
\kappa(x)\psi_{S_0}(yx)W'(a(y))\,dx\\
&=h(y)W'(a(y)),
\quad y\in F_{S_0}^{\times}.
\end{aligned}
\]
Consequently,
\[
W'_h(a(\gamma))=
\begin{cases}
1,&\gamma=\gamma_0,\\
0,&\gamma\in\Gamma\setminus\{\gamma_0\}
\end{cases}
\]

Since $\kappa=\bigotimes_{v\in S_0}\kappa_v$,
\[
W'_h=\bigotimes_{v\in S_0}W'_{v,h_v}.
\]
We may therefore choose a smooth vector 
\[
\xi_{S_0}^{(\gamma_0)}
=\bigotimes_{v\in S_0}\xi_v^{(\gamma_0)}
\]
such that
\[
R(g_{S_0})W_{\xi_{S_0}^{(\gamma_0)}}=W'_h.
\]
This proves the assertion.
\end{proof}

\begin{lemma}\label{lemma-nonvanishing-class}
Put $\mathfrak d=\operatorname{div}^{S_0}({\sf d})$.
For every admissible modulus $\mathfrak m$, there
exists an effective divisor $\mathfrak D$, supported
outside $S_0$, such that
\[
[\mathfrak D]=[\mathfrak n\mathfrak d^2]
\quad\text{in }\operatorname{Cl}(\o_F^{S_0}),
\qquad
(\mathfrak D,\mathfrak m)=1,
\qquad
\widetilde c^{S_0}(\mathfrak D\mathfrak d^{-1})
\ne0.
\]
\end{lemma}

\begin{proof}
Put $A=\o_F^{S_0}$ and
$T=\operatorname{Supp}(\mathfrak m)$.
Let $h$ be the order of the finite group
$\operatorname{Cl}(A)$.
For each $v\in T$, choose a sufficiently large
positive integer $e_v$ divisible by $h$, and a
character $\rho_v$ of $\o_v^\times$ of conductor
exponent $e_v$.
Choose a unitary id\`ele class character $\theta$
unramified outside $T\cup S_0$ such that
\[
\theta_v|_{\o_v^\times}=\rho_v
\qquad(v\in T),
\]
and put $\sigma=\pi\otimes\theta$.

Since $\pi_v$ is unramified for $v\in T$,
\[
L(s,\sigma_v)=L(s,\widetilde\sigma_v)=1,
\qquad
n(\sigma_v)=2e_v.
\]
Outside $T\cup S_0$, twisting by $\theta_v$
preserves the conductor.
Thus the conductor of $\sigma$ outside $S_0$ is
\[
\mathfrak N_\theta
=
\mathfrak n\prod_{v\in T}\mathfrak p_v^{2e_v}.
\]
By our choice of the $e_v$,
\begin{equation}\label{eqn-twist-conductor-class}
[\mathfrak N_\theta]=[\mathfrak n]
\quad\text{in }\operatorname{Cl}(A).
\end{equation}
Let ${\sf N}_\theta$ represent $\mathfrak N_\theta$,
with components $1$ at $S_0$.

For every id\`ele class character $\eta$ unramified
outside $S_0$, the character $\theta\eta$ is an
allowed twist in Theorem~\ref{mainthm}.
The original functional-equation scalar is
independent of $\eta$; denote it by
$\epsilon_\theta^{\mathrm{orig}}\ne0$.
Consequently, the base-case Mellin-inversion
argument applies to $\sigma$ and this family
of characters.
Only that argument is used here, not the
conclusion of Theorem~\ref{mainthm} for $\sigma$.

Choose normalized essential vectors for
$\sigma_v$ at $v\notin S_0$ and arbitrary smooth
vectors at $S_0$.
Denote the resulting functions by
$\phi^\sigma_{\xi_{S_0}}$ and
$\widetilde\phi^\sigma_{\xi_{S_0}}$, using the
original, unscaled dual function.
Applying the base-case identity at
$t={\sf d}^{-1}$ gives
\begin{equation}\label{eqn-twist-base}
\phi^\sigma_{\xi_{S_0}}(a({\sf d}^{-1}))
=
\epsilon_\theta^{\mathrm{orig}}\,
\widetilde\phi^\sigma_{\xi_{S_0}}
\bigl(a({\sf N}_\theta{\sf d})\bigr).
\end{equation}

Since ${\sf d}_{S_0}=1$, the Fourier expansion is
\[
\phi^\sigma_{\xi_{S_0}}(a({\sf d}^{-1}))
=
\sum_{\gamma\in A\setminus\{0\}}
c_\sigma^{S_0}(\gamma{\sf d}^{-1})
W^\sigma_{\xi_{S_0}}(a(\gamma)).
\]
Indeed, the support condition outside $S_0$
requires $\operatorname{ord}_v(\gamma)\geq0$.
The essential-vector normalization gives
\[
c_\sigma^{S_0}({\sf d}^{-1})=1.
\]
Apply Lemma~\ref{lemma-isolation} to $\sigma$,
with $\mathcal L=A$, $\gamma_0=1$, and
$g_{S_0}=1$.
We may choose a pure tensor $\xi_{S_0}$ so that
only the term $\gamma=1$ survives, with value $1$.
Hence
\[
\phi^\sigma_{\xi_{S_0}}(a({\sf d}^{-1}))=1.
\]
It follows from \eqref{eqn-twist-base} that
\[
\widetilde\phi^\sigma_{\xi_{S_0}}
\bigl(a({\sf N}_\theta{\sf d})\bigr)\ne0.
\]
Therefore there exists $\gamma\in F^\times$ such
that
\[
\widetilde c_\sigma^{S_0}
(\gamma{\sf N}_\theta{\sf d})\ne0.
\]

Put
\[
\mathfrak D
=
\operatorname{div}^{S_0}
(\gamma{\sf N}_\theta{\sf d}^{\,2}).
\]
By the support of the essential Whittaker
functions, this divisor is effective.
Equation~\eqref{eqn-twist-conductor-class} gives
\[
[\mathfrak D]
=
[\mathfrak N_\theta\mathfrak d^2]
=
[\mathfrak n\mathfrak d^2].
\]

For $v\in T$, the identity
$L(s,\widetilde\sigma_v)=1$ implies
\[
(L({\sf d}_v)W_{\widetilde\sigma_v})(a(y))
=
\mathbf 1_{\o_v^\times}(y).
\]
Equivalently,
\[
W_{\widetilde\sigma_v}(a(y))
=
\mathbf 1_{\o_v^\times}({\sf d}_v y).
\]
The local factor in $\widetilde c_\sigma^{S_0}$
is a nonzero scalar multiple of this function.
Thus its nonvanishing implies
\[
\operatorname{ord}_v
(\gamma{\sf N}_{\theta,v}{\sf d}_v^2)=0
\qquad(v\in T).
\]
Consequently, $(\mathfrak D,\mathfrak m)=1$.

Outside $T\cup S_0$, the normalized essential
diagonal coefficients for
$\widetilde\sigma_v
=\widetilde\pi_v\otimes\theta_v^{-1}$
differ from those for $\widetilde\pi_v$ by
nonzero character factors.
At $v\in T$, the argument
$y=\gamma{\sf N}_{\theta,v}{\sf d}_v$ satisfies
${\sf d}_v y\in\o_v^\times$.
At such an argument, both normalized essential
Whittaker functions have value $1$.
The corresponding dual local factors are
therefore nonzero.
It follows that
\[
\widetilde c^{S_0}
(\gamma{\sf N}_\theta{\sf d})\ne0.
\]
Since
\[
\operatorname{div}^{S_0}
(\gamma{\sf N}_\theta{\sf d})
=
\mathfrak D\mathfrak d^{-1},
\]
the assertion follows.
\end{proof}

\subsection{Conclusion of the proof of the main theorem}

We keep $\mathfrak H$ and its
representing id\`ele ${\sf h}$ as in \S\ref{ss-aut}.
Fix an admissible modulus $\mathfrak m$ coprime
to $\mathfrak H$ such that
$[\mathfrak m]_{\mathfrak n\mathfrak H}=1$ in
$\operatorname{Cl}_{\mathfrak n\mathfrak H}(A)$.
Thus there exists $M\in F^\times$ satisfying
\[
\operatorname{div}^{S_0}(M)=\mathfrak m,
\qquad
M\equiv1\pmod{\mathfrak n\mathfrak H}.
\]
We retain ${\sf f}_1,{\sf e}_1,{\sf e}'_1,t_1$
from \S\ref{ss-aut}, and put
\[
z_{\mathfrak m}={\sf n}t_1^{-1}{\sf e}'_1,
\qquad
g_v=a(t_{1,v}^{-1}{\sf f}_{1,v})k_v
\quad(v\in S_0).
\]
Write $g_{S_0}=(g_v)_{v\in S_0}$ and set
\[
\Gamma_{\mathfrak m}
=
\left\{
\gamma\in F^\times:
\operatorname{ord}_v
({\sf d}_v\gamma{\sf n}_v{\sf h}_v
 {\sf f}_{\mathfrak m,v})
\geq0
\text{ for every }v\notin S_0
\right\}.
\]
We suppress the dependence of
$\Gamma_{\mathfrak m}$ on the fixed divisor
$\mathfrak H$.
The coefficient
$\widetilde c^{S_0}
(\gamma{\sf n}{\sf h}{\sf f}_{\mathfrak m})$
vanishes outside $\Gamma_{\mathfrak m}$.
For $\gamma\in\Gamma_{\mathfrak m}$, the argument
of $S_{\mathfrak m}$ is the image of
\[
({\sf d}_v\gamma{\sf n}_v{\sf h}_v
 {\sf f}_{\mathfrak m,v})_{v\in S}
\]
in $R_{\mathfrak m}$.

In \eqref{eqn-R}, take $b=1$ and replace $t$ by
$t{\sf f}_{\mathfrak m}$. At every $v\notin S_0$,
choose
\[
{\sf f}_v=1,\qquad {\sf e}'_v=0,\qquad
{\sf e}_v=
\begin{cases}
0,&v\in S,\\
{\sf h}_v,&v\notin S.
\end{cases}
\]
At the places in $S_0$, take
\[
({\sf f}_v,{\sf e}_v,{\sf e}'_v)
=
({\sf f}_{1,v},{\sf e}_{1,v},{\sf e}'_{1,v}).
\]
Finally, set $t=t_1$.

With these choices, the matrix on the left is
\[
\begin{pmatrix}
t_1{\sf f}_1&t_1{\sf e}_1\\
0&1
\end{pmatrix},
\]
and the translated matrix on the dual side is
\[
\begin{pmatrix}
{\sf n}t_1^{-1}{\sf f}_1&
{\sf n}t_1^{-1}{\sf e}'_1\\
0&1
\end{pmatrix}.
\]
Consequently,
\begin{align}
&\phi_{\xi_{S_0}}\!\left(
\begin{psmallmatrix}
t_1{\sf f}_1&t_1{\sf e}_1\\&1
\end{psmallmatrix}\right)
\notag\\
&\quad=
\chi_\pi({\sf f}_1)
\sum_{\gamma\in\Gamma_{\mathfrak m}}
\widetilde c^{S_0}
(\gamma{\sf n}{\sf h}{\sf f}_{\mathfrak m})
S_{\mathfrak m}
({\sf d}\gamma{\sf n}{\sf h}{\sf f}_{\mathfrak m})
\psi(\gamma z_{\mathfrak m})
\widetilde W_{\xi_{S_0}}
(a(\gamma)g_{S_0}).
\label{eq:s43-weighted}
\end{align}

On the other hand, expanding the right-hand side
of Proposition~\ref{prop-special} gives the same
expression without the factor $S_{\mathfrak m}$.
Subtracting these identities and dividing by the
nonzero scalar $\chi_\pi({\sf f}_1)$ yields
\begin{equation}\label{eq:s43-comparison}
\begin{aligned}
\sum_{\gamma\in\Gamma_{\mathfrak m}}
&\widetilde c^{S_0}
(\gamma{\sf n}{\sf h}{\sf f}_{\mathfrak m})
\left(
S_{\mathfrak m}
({\sf d}\gamma{\sf n}{\sf h}{\sf f}_{\mathfrak m})-1
\right)\\
&\qquad\times
\psi(\gamma z_{\mathfrak m})
\widetilde W_{\xi_{S_0}}
(a(\gamma)g_{S_0})
=0.
\end{aligned}
\end{equation}
This holds for every smooth vector
$\xi_{S_0}$ in the space of $\pi_{S_0}$.

The set $\Gamma_{\mathfrak m}$ is the set of
nonzero elements of the fractional ideal
\[
\mathcal L_{\mathfrak m}
=
(\mathfrak d\mathfrak n\mathfrak H\mathfrak m)^{-1}
\]
of $A$.
Apply Lemma~\ref{lemma-isolation} to
$\widetilde\pi$, with this fractional ideal,
$\gamma_0\in\Gamma_{\mathfrak m}$ and $g_{S_0}$.
The resulting pure tensor of Whittaker functions
can be realized as $\widetilde W_{\xi_{S_0}}$
for a suitable pure tensor $\xi_{S_0}$.
Applying \eqref{eq:s43-comparison} to this vector,
we obtain
\[
\widetilde c^{S_0}
(\gamma_0{\sf n}{\sf h}{\sf f}_{\mathfrak m})
\left(
S_{\mathfrak m}
({\sf d}\gamma_0{\sf n}{\sf h}{\sf f}_{\mathfrak m})-1
\right)
\psi(\gamma_0z_{\mathfrak m})=0.
\]
Since the additive character never vanishes and
$\gamma_0$ was arbitrary, it follows that
\begin{equation}\label{eq:s43-isolated}
\widetilde c^{S_0}
(\gamma{\sf n}{\sf h}{\sf f}_{\mathfrak m})
\left(
S_{\mathfrak m}
({\sf d}\gamma{\sf n}{\sf h}{\sf f}_{\mathfrak m})-1
\right)=0
\qquad(\gamma\in\Gamma_{\mathfrak m}).
\end{equation}

\begin{lemma}\label{lemma-S-compatibility}
For any admissible moduli $\mathfrak m\mid\m'$, suppose 
\[
\tau={\sf f}_{\m'}{\sf f}_{\mathfrak m}^{-1}.
\]
If $y_v\in\o_v$ for every $v\mid\m'$, then
\[
S_{\m'}(\tau y)=S_{\mathfrak m}(y),
\]
where the arguments are interpreted in
$R_{\m'}$ and $R_{\mathfrak m}$, respectively.
\end{lemma}

\begin{proof}
Let
\[
r:\mathcal K_{\mathfrak m'}\longrightarrow
  \mathcal K_{\mathfrak m}
\]
be the natural reduction map, and put $J=\ker r$.
This map is surjective, so that
$N_{\mathfrak m'}=|J|N_{\mathfrak m}$. First, note that
\[
\Psi_{\mathfrak m'}((1-u')\tau y)
=
\Psi_{\mathfrak m}((1-r(u'))y).
\]
Indeed, if $a_v=\operatorname{ord}_v(\mathfrak m)$
and $a'_v=\operatorname{ord}_v(\mathfrak m')$, then $\tau_v=\varpi_v^{a'_v-a_v}$.
Thus, at $v\mid\mathfrak m$,
\[
\begin{aligned}
&\psi_v\bigl(
{\sf d}_v^{-1}\varpi_v^{-a'_v}
(1-u'_v)\tau_v y_v
\bigr)\\
&\qquad=
\psi_v\bigl(
{\sf d}_v^{-1}\varpi_v^{-a_v}
(1-u'_v)y_v
\bigr).
\end{aligned}
\]
At $v\mid\mathfrak m'$ with $v\nmid\mathfrak m$,
the corresponding factor is
\[
\psi_v\bigl({\sf d}_v^{-1}(1-u'_v)y_v\bigr)=1,
\]
since $(1-u'_v)y_v\in\o_v$.
Thus 
\begin{align*}
S_{\mathfrak m'}(\tau y)
&=\frac{1}{N_{\mathfrak m'}}
  \sum_{u'\in\mathcal K_{\mathfrak m'}}
  \widehat\epsilon_{\mathfrak m'}(u')
  \Psi_{\mathfrak m'}((1-u')\tau y)\\
&=\frac{1}{N_{\mathfrak m'}}
  \sum_{u\in\mathcal K_{\mathfrak m}}
  \Psi_{\mathfrak m}((1-u)y)
  \sum_{\{u':r(u')=u\}}\widehat\epsilon_{\mathfrak m'}(u').
\end{align*}
To evaluate the inner sum, expand
$\widehat\epsilon_{\mathfrak m'}$:
\[
\sum_{\{u':r(u')=u\}}\widehat\epsilon_{\mathfrak m'}(u')
=
\sum_{\rho'\in\widehat{\mathcal K}_{\mathfrak m'}}
\epsilon_{[\rho']}
\sum_{\{u':r(u')=u\}}\rho'^{-1}(u').
\]
The inner sum vanishes
unless $\rho'$ is trivial on $J$. In that case,
$\rho'=\rho\circ r$ for a unique
$\rho\in\widehat{\mathcal K}_{\mathfrak m}$, and the
inner sum equals $|J|\rho^{-1}(u)$.
Moreover, since there is a common lift $\omega\in X(\m)$ of $\rho$ and $\rho\circ r$, we have 
\[
\epsilon_{[\rho\circ r]}
=\epsilon_{[\rho]}.
\]
Consequently,
\[
\sum_{\{u':r(u')=u\}}\widehat\epsilon_{\mathfrak m'}(u')
=
|J|\widehat\epsilon_{\mathfrak m}(u).
\]
Substitution now gives
\begin{align*}
S_{\mathfrak m'}(\tau y)
&=\frac{|J|}{N_{\mathfrak m'}}
  \sum_{u\in\mathcal K_{\mathfrak m}}
  \widehat\epsilon_{\mathfrak m}(u)
  \Psi_{\mathfrak m}((1-u)y)\\
&=\frac{1}{N_{\mathfrak m}}
  \sum_{u\in\mathcal K_{\mathfrak m}}
  \widehat\epsilon_{\mathfrak m}(u)
  \Psi_{\mathfrak m}((1-u)y)
=S_{\mathfrak m}(y),
\end{align*}
as required.
\end{proof}

\begin{proposition}\label{prop-S-one}
For every admissible modulus $\mathfrak m$,
\[
S_{\mathfrak m}(x)=1
\qquad \forall x\in R_{\mathfrak m}.
\]
\end{proposition}

\begin{proof}
Fix $x\in R_{\m}$. Choose
$\mathfrak D$ as in Lemma~\ref{lemma-nonvanishing-class}.
Thus
\[
[\mathfrak D]=[\mathfrak n\mathfrak d^2],
\qquad
(\mathfrak D,\mathfrak m)=1,
\qquad
\widetilde c^{S_0}(\mathfrak D\mathfrak d^{-1})\ne0.
\]

Choose an effective divisor $\mathfrak H$ such that
\[
[\mathfrak H]=[\mathfrak d]
\quad\text{in }\operatorname{Cl}(A),
\qquad
(\mathfrak H,\mathfrak m\mathfrak n\mathfrak D)=1.
\]
Such a choice is possible by the Chinese remainder
theorem: choose $0\ne b\in\mathfrak d^{-1}$ whose
image in
$\mathfrak d^{-1}/\mathfrak p\mathfrak d^{-1}$
is nonzero at every prime $\mathfrak p$ dividing
$\mathfrak m\mathfrak n\mathfrak D$, and put
$\mathfrak H=(b)\mathfrak d$.
Let ${\sf h}$ represent $\mathfrak H$ as in \S\ref{ss-aut}.

Since
$[\mathfrak d\mathfrak n\mathfrak H]=[\mathfrak D]$,
choose $\delta\in F^\times$ such that
\[
(\delta)\mathfrak d\mathfrak n\mathfrak H
=
\mathfrak D
\quad\text{as divisors of }A.
\]
Another application of the Chinese remainder theorem to the
$A$-module $\mathfrak m^{-1}$, choose
$0\ne a\in\mathfrak m^{-1}$ such that
\[
{\sf d}_v\delta a\,{\sf n}_v{\sf h}_v
{\sf f}_{\mathfrak m,v}
\equiv x_v
\pmod{\mathfrak p_v^{\operatorname{ord}_v(\mathfrak m)}}
\qquad(v\mid\mathfrak m),
\]
\[
a\equiv1
\pmod{\mathfrak p_v^{
\operatorname{ord}_v(\mathfrak n\mathfrak H)}}
\qquad(v\mid\mathfrak n\mathfrak H),
\]
and
\[
a\equiv1\pmod{\mathfrak p_v}
\qquad
(v\mid\mathfrak D,\;v\nmid\mathfrak n\mathfrak H).
\]
These conditions are compatible: the supports of
$\mathfrak m$ and $\mathfrak n\mathfrak H\mathfrak D$
are disjoint, and
${\sf d}_v\delta{\sf n}_v{\sf h}_v$ is a unit
for every $v\mid\mathfrak m$, since
\[
(\delta)\mathfrak d\mathfrak n\mathfrak H
=\mathfrak D.
\]
Thus, the first family prescribes residue classes
in the localizations of $\mathfrak m^{-1}$.

In particular, $a$ is integral at every place
outside $S_0\cup\operatorname{Supp}(\mathfrak m)$
and is a unit at every prime dividing
$\mathfrak n\mathfrak H\mathfrak D$. Hence
\[
\mathfrak I=(a)\mathfrak m
\]
is integral and satisfies
\[
(\mathfrak I,\mathfrak n\mathfrak H\mathfrak D)=1.
\]
Moreover, since $a\equiv1\pmod{\mathfrak n\mathfrak H}$,
\[
[\mathfrak I]_{\mathfrak n\mathfrak H}
=
[\mathfrak m]_{\mathfrak n\mathfrak H}.
\]
Putting
\[
\gamma=\delta a,
\qquad
y={\sf d}\gamma{\sf n}{\sf h}{\sf f}_{\mathfrak m},
\]
we have
\[
\operatorname{div}^{S_0}(y)=\mathfrak D\mathfrak I.
\]
Consequently, $y$ is integral at every place outside
$S_0$ and satisfies
\[
y_v\equiv x_v\pmod{\mathfrak p_v^{n_v}}
\qquad(v\mid\mathfrak m),
\]
where $n_v=\operatorname{ord}_v(\mathfrak m)$.

Let $d$ be the order of
$[\mathfrak m]_{\mathfrak n\mathfrak H}$ in
$\operatorname{Cl}_{\mathfrak n\mathfrak H}(A)$.
Choose an integer $E\geq2$ divisible by $d$.
This ensures
\[
[\mathfrak m]_{\mathfrak n\mathfrak H}^{E}=1.
\]
Since
\[
[\mathfrak I]_{\mathfrak n\mathfrak H}
=
[\mathfrak m]_{\mathfrak n\mathfrak H},
\]
it follows that
\[
\mathfrak M=\mathfrak m\mathfrak I^{E-1}
\]
is trivial in
$\operatorname{Cl}_{\mathfrak n\mathfrak H}(A)$.
Consequently, there exists $M\in F^\times$ such that
\[
\operatorname{div}^{S_0}(M)=\mathfrak M,
\qquad
M\equiv1\pmod{\mathfrak n\mathfrak H},
\]
as required by Proposition~\ref{prop-special}.
Also,
$(\mathfrak M,\mathfrak n\mathfrak H)=1$.

Second, we also require
\[
(L({\sf d}_v)W_{\widetilde\pi_v})
\left(a(\varpi_v^{Er_v})\right)\ne0
\qquad(v\mid\mathfrak I),
\]
where $r_v=\operatorname{ord}_v(\mathfrak I)$.
Since $(\mathfrak I,\mathfrak n)=1$, the
representation $\pi_v$ is unramified.
Let $\alpha_v,\beta_v$ be its Satake parameters,
and write $\zeta_v=\alpha_v/\beta_v$.
If $\alpha_v\ne\beta_v$, then
\[
\lambda_{Er_v}
=
\frac{\alpha_v^{Er_v+1}-\beta_v^{Er_v+1}}
     {\alpha_v-\beta_v}
=
\beta_v^{Er_v}
\frac{\zeta_v^{Er_v+1}-1}{\zeta_v-1}.
\]
Thus $\lambda_{Er_v}=0$ precisely when
$\zeta_v^{Er_v+1}=1$.
If $\zeta_v$ is not a root of unity, this cannot
occur. If $\zeta_v\ne1$ has finite order $d_v$,
we additionally require $d_v\mid E$. Then
\[
\zeta_v^{Er_v+1}=\zeta_v\ne1,
\]
so again $\lambda_{Er_v}\ne0$.
Finally, if $\alpha_v=\beta_v$, then
\[
\lambda_{Er_v}
=(Er_v+1)\alpha_v^{Er_v}\ne0.
\]

There are only finitely many primes dividing
$\mathfrak I$. We may therefore choose $E\geq2$
to be divisible by $d$ and by every order $d_v$
arising above. This choice ensures both that
\[
[\mathfrak M]_{\mathfrak n\mathfrak H}=1
\]
and that
\[
\lambda_{Er_v}\ne0
\qquad(v\mid\mathfrak I).
\]
The corresponding coefficients for
$\widetilde\pi_v$ are also nonzero, since
\[
\widetilde\lambda_{Er_v}
=
(\alpha_v\beta_v)^{-Er_v}\lambda_{Er_v}.
\]

We can now use \eqref{eq:s43-isolated} with
$\mathfrak m=\mathfrak M$ and the auxiliary
divisor $\mathfrak H$ chosen above.
Since
\[
\operatorname{div}^{S_0}
({\sf d}\gamma{\sf n}{\sf h}{\sf f}_{\mathfrak M})
=
\mathfrak D\mathfrak I^E,
\]
we have $\gamma\in\Gamma_{\mathfrak M}$ and
\[
\operatorname{div}^{S_0}
(\gamma{\sf n}{\sf h}{\sf f}_{\mathfrak M})
=
\mathfrak D\mathfrak I^E\mathfrak d^{-1}.
\]
We compare the local factors of
\[
\widetilde c^{S_0}
(\mathfrak D\mathfrak I^E\mathfrak d^{-1})
\]
with those of
$\widetilde c^{S_0}(\mathfrak D\mathfrak d^{-1})\ne0$.
Indeed, for $v\notin S_0$ with $v\nmid\mathfrak I$,
we have
\[
\operatorname{ord}_v
(\mathfrak D\mathfrak I^E\mathfrak d^{-1})
=
\operatorname{ord}_v(\mathfrak D\mathfrak d^{-1}).
\]
Since the local factor depends only on this
valuation, it is unchanged and remains nonzero.

For $v\mid\mathfrak I$, the coprimality
$(\mathfrak D,\mathfrak I)=1$ gives
\[
\operatorname{ord}_v(\mathfrak D)=0,
\qquad
\operatorname{ord}_v(\mathfrak D\mathfrak I^E)=Er_v.
\]
Thus, after the translation by $L({\sf d}_v)$,
the relevant essential Whittaker value is at
$a(\varpi_v^{Er_v})$.
Our choice of $E$ ensures that this value,
and hence the corresponding dual local factor,
is nonzero. Consequently,
\[
\widetilde c^{S_0}
(\gamma{\sf n}{\sf h}{\sf f}_{\mathfrak M})
=
\widetilde c^{S_0}
(\mathfrak D\mathfrak I^E\mathfrak d^{-1})
\ne0.
\]

From \eqref{eq:s43-isolated} with modulus
$\mathfrak M$, we have
\[
\widetilde c^{S_0}
(\gamma{\sf n}{\sf h}{\sf f}_{\mathfrak M})
\left(
S_{\mathfrak M}
({\sf d}\gamma{\sf n}{\sf h}{\sf f}_{\mathfrak M})-1
\right)=0,
\]
from which we conclude
\[
S_{\mathfrak M}
({\sf d}\gamma{\sf n}{\sf h}{\sf f}_{\mathfrak M})=1.
\]

Finally, put
\[
\tau={\sf f}_{\mathfrak M}{\sf f}_{\mathfrak m}^{-1},
\qquad
y={\sf d}\gamma{\sf n}{\sf h}{\sf f}_{\mathfrak m}.
\]
Since
$\operatorname{div}^{S_0}(y)=\mathfrak D\mathfrak I$,
the integrality requirement in
Lemma~\ref{lemma-S-compatibility} is satisfied.
That lemma and the choice of $a$ give
\[
1=S_{\mathfrak M}(\tau y)
 =S_{\mathfrak m}(y)
 =S_{\mathfrak m}(x).
\]
Since $x$ was arbitrary, the proposition follows.
\end{proof}
\begin{corollary}\label{cor-eps}
For every admissible modulus $\mathfrak m$ and
$\omega\in X(\mathfrak m)$, the normalized scalar
$\epsilon_{[\omega]}$ equals $1$.
\end{corollary}
\begin{proof}
By the proposition, $S_{\mathfrak m}(x)=1$ for every
$x\in R_{\mathfrak m}$. Fix $u_0\in\mathcal K_{\mathfrak m}$.
Multiplying the defining identity for $S_{\mathfrak m}(x)$
by $\Psi_{\mathfrak m}((u_0-1)x)$ and averaging over
$R_{\mathfrak m}$, we obtain
\begin{align*}
\mathbf 1_{\{1\}}(u_0)
&=\frac{1}{|R_{\mathfrak m}|}
  \sum_{x\in R_{\mathfrak m}}
  \Psi_{\mathfrak m}((u_0-1)x)\\
&=\frac{1}{N_{\mathfrak m}}
  \sum_{u\in\mathcal K_{\mathfrak m}}
  \widehat\epsilon_{\mathfrak m}(u)
  \left(
  \frac{1}{|R_{\mathfrak m}|}
  \sum_{x\in R_{\mathfrak m}}
  \Psi_{\mathfrak m}((u_0-u)x)
  \right)\\
&=\frac{\widehat\epsilon_{\mathfrak m}(u_0)}
        {N_{\mathfrak m}},
\end{align*}
where the first and last equalities follow from
additive character orthogonality and the non-degeneracy
of the pairing $(a,x)\mapsto\Psi_{\mathfrak m}(ax)$.
Thus
\[
\widehat\epsilon_{\mathfrak m}(u)
=N_{\mathfrak m}\mathbf 1_{\{1\}}(u).
\]
From this, using multiplicative Fourier inversion, we obtain
\[
\epsilon_{[\omega]}
=
\frac{1}{N_{\mathfrak m}}
\sum_{u\in\mathcal K_{\mathfrak m}}
\widehat\epsilon_{\mathfrak m}(u)\omega_S(u)
=1.
\]
\end{proof}

We use the following variant of the converse theorem
of \cite{CKPSS2} for $\mathrm{GL}(2)$, allowing a fixed
nonzero scalar in the functional equations.
The same modification applies to $\mathrm{GL}(N)$,
provided the scalar is independent of the twisting
representation and the remaining hypotheses of the
cited converse theorem are satisfied.

\begin{lemma}\label{lemma-constant-converse}
Let $\pi=\bigotimes_v'\pi_v$ be an irreducible admissible
generic representation of ${\rm GL}_2(\A_F)$ with
automorphic central character, and let $\Sigma$ be a
finite set of finite places.
Suppose that, for every unitary idele class character
$\chi$ unramified at all places in $\Sigma$, the completed
$L$-functions
\[
\Lambda(s,\pi\times\chi),
\qquad
\Lambda(s,\widetilde\pi\times\chi^{-1})
\]
converge absolutely in a right half-plane, extend to
entire functions bounded in finite vertical strips,
and satisfy
\[
\Lambda(s,\pi\times\chi)
=
c\,\epsilon(s,\pi\times\chi,\psi)
\Lambda(1-s,\widetilde\pi\times\chi^{-1}),
\]
where $c\ne0$ is independent of $\chi$.
Then there exists a unique isobaric automorphic
representation $\Pi$ of $\GL_2(\A_F)$ such that
$\Pi_v\simeq\pi_v$ at every nonarchimedean place
$v\notin\Sigma$ where $\pi_v$ is unramified.
\end{lemma}

\begin{proof}
We follow the proof of the converse theorem in
\cite[\S~2]{CKPSS2}, with $N=2$ and auxiliary character
$\eta=1$, retaining the scalar in the functional
equations. We use the notation and local test vectors of that proof.
In particular, fix the prescribed vectors
$\xi_\Sigma^0$ at the places in $\Sigma$, allow
$\xi^\Sigma$ to vary, and put
$\xi=\xi_\Sigma^0\otimes\xi^\Sigma$.
Let $U_\xi$ and $V_\xi$ be the two mirabolic
Whittaker sums constructed there, and let
\[
G_0=K_{0,\Sigma}(\mathfrak a)\,G^\Sigma
\]
be the subgroup occurring in their proof. The functional equations and Mellin inversion now
give
\[
U_\xi(g)=c\,V_\xi(g)\qquad(g\in G_0),
\]
in place of $U_\xi(g)=V_\xi(g)$.
The same scalar occurs for every $\xi^\Sigma$.

Let
\[
P_0=P(F)\cap G_0,\qquad
Q_0=Q(F)\cap G_0,
\]
where $P$ and $Q$ are the mirabolic subgroups used
in the cited proof.
The function $U_\xi$ is left $P_0$-invariant by
construction. For $q\in Q_0$ and $g\in G_0$,
the left $Q(F)$-invariance of $V_\xi$ gives
\[
U_\xi(qg)=c\,V_\xi(qg)=c\,V_\xi(g)=U_\xi(g).
\]
Thus $U_\xi|_{G_0}$ is invariant under $P_0$ and $Q_0$.
It is also invariant under rational central elements,
since the central character is automorphic.
The remaining generation, approximation, and
representation-theoretic arguments in the cited proof
are unchanged by the fixed scalar.
Thus we obtain an irreducible automorphic representation
$\Pi'$ with $\Pi'_v\simeq\pi_v$ for every $v\notin\Sigma$.

Finally, let $\Pi$ be the isobaric automorphic representation
associated to the cuspidal support of $\Pi'$.
By the argument in the proof of
\cite[Lemma~3.5]{BK2}, this replacement
preserves the unramified local components. Hence $\Pi_v\simeq\Pi'_v\simeq\pi_v$ at every
nonarchimedean place $v\notin\Sigma$ where
$\pi_v$ is unramified. The uniqueness of $\Pi$ follows from strong multiplicity
one for isobaric automorphic representations.
\end{proof}

By Corollary~\ref{cor-eps}, the normalized scalar factors 
are all equal to~$1$; see Remark~\ref{rmk-normalization}
in the Introduction. Undoing the normalization, we find
that the original factors $\epsilon_{[\chi]}$ are
equal to a fixed nonzero constant, independent of $\chi$.
The required boundedness in vertical strips
follows by Phragm\'en--Lindel\"of.
Taking
\[
\Sigma=\{v<\infty:\pi_v\text{ is ramified}\},
\]
the hypotheses of Lemma~\ref{lemma-constant-converse}
are therefore satisfied, and its conclusion gives
the required isobaric automorphic representation.
This completes the proof of Theorem~\ref{mainthm}.

\bibliographystyle{amsplain}

\providecommand{\bysame}{\leavevmode\hbox to3em{\hrulefill}\thinspace}
\providecommand{\MR}{\relax\ifhmode\unskip\space\fi MR }
\providecommand{\MRhref}[2]{%
  \href{http://www.ams.org/mathscinet-getitem?mr=#1}{#2}%
}

\end{document}